\documentclass[12pt,reqno]{amsart}
\usepackage[dvips]{graphics}

\newtheorem{theorem}{{\bf Theorem}}[section]
\newtheorem*{maintheorem}{{\bf Main Theorem}}
\newtheorem{proposition}[theorem]{{\bf Proposition}}
\newtheorem*{proposition*}{{\bf Proposition}}
\newtheorem{definition}[theorem]{{\bf Definition}}
\newtheorem{lemma}[theorem]{{\bf Lemma}}
\newtheorem{lemma*}{{\bf Lemma}}
\newtheorem{example}[theorem]{{\bf Example}}

\begin{document}

%
%

\title[Circle Packings On Surfaces]{%
	CIRCLE PACKINGS ON SURFACES WITH PROJECTIVE STRUCTURES}

\author[S. Kojima]{%
	Sadayoshi Kojima} 
\address{%
        Department of Mathematical and Computing Sciences \\
        Tokyo Institute of Technology \\
        Ohokayama, Meguro \\    
        Tokyo 152-8552 Japan} 
\email{%
        sadayosi@is.titech.ac.jp}  

\author[S. Mizushima]{%
	Shigeru Mizushima} 
\address{%
        Department of Mathematical and Computing Sciences \\
        Tokyo Institute of Technology \\
        Ohokayama, Meguro \\    
        Tokyo 152-8552 Japan} 
\email{%
        mizusima@is.titech.ac.jp}  

\author[S. P. Tan ]{%
	Ser Peow Tan} 
\address{%
        Department of Mathematics \\
        National University of Singapore \\
        Singapore 117543, \\    
        Singapore} 
\email{%
        mattansp@nus.edu.sg}  


\thanks{The third author gratefully acknowledges support of the National University of Singapore academic research grant R-146-000-031-112 and  the Japan Society for the Promotion of Sciences for support for a visit to the Tokyo Institute of Technology in May 2000.}

%
%

\begin{abstract}
The Andreev-Thurston theorem states that 
for any triangulation of a closed orientable 
surface $\Sigma_g$ of genus $g$ which is covered by a simple graph 
in the universal cover, 
there exists a unique metric of curvature $1, 0$ or $-1$ on the surface 
depending on whether $g=0, 1$ or $ \ge 2$ such that 
the surface with this metric admits 
a circle packing with combinatorics given by the triangulation.  
Furthermore, the circle packing is essentially rigid, 
that is, unique up to conformal automorphisms of the surface 
isotopic to the identity.  
In this paper, we consider projective 
structures on the surface $\Sigma_g$  
where circle packings are also defined.  
We show that the space of projective structures 
on a surface of genus $g \ge 2$ which admits 
a circle packing by one circle is homeomorphic 
to  ${\mathbb R}^{6g-6}$  and furthermore that  
the circle packing is rigid on such surfaces.
\end{abstract}

\maketitle

%
%

\section{Introduction}

This paper is motivated by the interplay between 
Mostow rigidity and hyperbolic Dehn surgery theory 
and is concerned with rigidity properties as well as 
deformation theory on the space of projective structures on surfaces.  
Mostow rigidity concerns the rigidity property of say,  
complete hyperbolic structures on manifolds of dimension $\geq 3$ based 
solely on the topological data of the manifold. On the other hand,  
if we relax the completeness condition in dimension three in particular, 
there is a non-trivial deformation theory.  
Hyperbolic Dehn surgery theory basically tells us how nice 
deformations are parameterized.  

Mostow rigidity fails in dimension two and 
the deformation space of hyperbolic structures on 
a closed orientable surface $\Sigma_g$ of genus $g \geq 2$ is 
homeomorphic to ${\mathbb R}^{6g-6}$.  
However the addition of certain ``tight'' topological data 
corresponding to circle packings can impose some rigidity properties 
as follows. 

Roughly, 
a circle packing on  $\Sigma_g$  is 
a collection of closed disks on the surface such that 
the interiors of any two distinct disks are disjoint and 
the complement of the disks in $\Sigma_g$ consists of disjoint, 
triangular interstices. 
One obtains a cellular decomposition of $\Sigma_g$ by assigning 
a vertex for each disk, 
an edge joining two vertices for each point of 
tangency of two (not necessarily distinct) disks and 
a triangle for each triangular interstice.  
We call the graph on $\Sigma_g$ consisting of the vertices and 
edges obtained in this way  the nerve of the circle packing.  
It  triangulates  $\Sigma_g$  in a generalized sense. 

The preimage of the nerve in the universal cover   
contains no loops by one or two edges, 
because otherwise, 
one would find a self-tangent circle or two circles 
with two points of tangency on the plane, 
which both would be absurd. 
In other words, 
the nerve must be covered by a simple graph 
in the graph theoretic sense.  

Now, the Andreev-Thurston theorem states that the surface  $\Sigma_g$  
together with a graph $\tau$ with the property above 
determines a unique metric of 
constant curvature on  $\Sigma_g$  such that the surface 
with this metric admits a circle packing whose nerve is isotopic to  $\tau$.  
Furthermore, this circle packing is rigid, that is, it is unique  up to conformal automorphisms of the surface isotopic to the identity.  
This result  was proved by Andreev in \cite{And} for the sphere and by Thurston in \cite{Thu} (see also \cite{BeSte1} and \cite{Col} for other proofs) for higher genus.  
Such a rigidity property implies that the set of points in the Teichm\"uller space which can be circle-packed is a countable set (Brooks proved in \cite{Brks} that this set is dense).  
We are interested in finding some relaxation of the conditions so that non-trivial deformations can occur.
 
A priori, the notion of a circle on a surface appears to be a metric notion, that is, a circle is the set of points equidistant from the center,  and  metrics of constant curvature form a very natural context for the study of circle packings. 
However there is a more general class of  geometric structures on surfaces called the complex projective structures in which the notion of a circle still makes sense.  
A complex projective structure (henceforth abbreviated to projective structure in this paper) is a $(G,X)$ structure (see \cite{Thu}) with the group $G=PSL_2({\mathbb C})$ and the model space $X=\hat{\mathbb C}$, the Riemann sphere.  
Two structures $(\Sigma_g, \mu_1)$ and $(\Sigma_g, \mu_2)$ are equivalent if there is a projective isomorphism isotopic to the identity between the two.  
The space of projective structures up to equivalence is denoted by ${\mathcal P}_g$  and is known to be homeomorphic to ${\mathbb R}^{12g-12}$ ($g \geq 2$).  
The uniformization defines a bundle projection
\begin{equation*} 
	u : {\mathcal P}_g \to {\mathcal T}_g
\end{equation*}
to the Teichm\"uller space ${\mathcal T}_g$, 
the space of all conformal structures on  $\Sigma_g$.  
The uniformization admits a canonical section  
\begin{equation*} 
	s : {\mathcal T}_g \to {\mathcal P}_g
\end{equation*}
which assigns to each conformal structure 
the equivalent hyperbolic structure, 
see for example \cite{Kmsm}. 
Since any projective transformation $\gamma \in PSL_2({\mathbb C})$ maps a circle on the Riemann sphere to a circle, the circle (and disk) makes sense as a geometric object on a surface with a projective structure.  Moreover, it is obvious that the nerve of a packing in a surface with projective structure again must be covered by a simple graph in the universal cover.     
The two main questions we are interested in are:

\smallskip

\noindent {\bf Question 1}. {\it For a fixed triangulation $\tau$ on  $\Sigma_g$  which is covered by a simple graph in the universal cover, what is the deformation space of projective structures which admit a circle packing with nerve isotopic to $\tau$?  
}

\medskip

\noindent {\bf Question 2}. {\it For a surface with a fixed projective structure which admits a circle packing with nerve $\tau$, is the circle packing rigid?
}
\medskip

Note that the Andreev-Thurston theorem implies that the intersection of this deformation space with  $s({\mathcal T}_g)$  is a point. 
The main results in this paper is the complete answer to the two questions above in the case where $\tau$ has only one vertex, that is, the circle packing of the surface consists of only one circle.

Our approach differs considerably from that used by most other authors working in the field of circle packings.  
In most previous papers, the authors use the radii of the circles as the parameters of the vertices of $\tau$.  
In our case, since a projective structure is not a metric structure, the radius does not make sense, neither does the center of the circle. Indeed, a careful analysis of the problem we are interested in shows that many of the methods and basic lemmas of circle packings on hyperbolic (or euclidean) surfaces do not apply here since in general, our structures have very complicated developing maps which are not covering maps onto the image and the holonomy group is in general not discrete.

To attack our problem, 
we introduce a cross ratio type invariant, 
which is a positive real number, 
for a certain configuration of 4 circles. 
It determines a configuration 
up to multiplication by elements 
of $PSL_2({\mathbb C})$. 
Suppose that we have a circle packing  $P$  on a surface with 
projective structure  $(\Sigma_g, \mu)$  with nerve  $\tau$.  
Then we define a map  
\begin{equation*} 
	{\bf c}_P : E(\tau) \to {\mathbb R}_+ 
\end{equation*} 
by assigning to each edge the cross ratio type invariant of its neighboring 4 circles in the universal cover, 
where  $E(\tau)$  is the set of edges of  $\tau$.  
We call  ${\bf c}_P$  a cross ratio parameter of the packing  $P$  on 
a surface with projective structure  $(\Sigma_g, \mu)$.  
It will be clear by defining the above 
more precisely in \S 2 that 
the cross ratio parameter is a complete projective 
invariant of the pair of a projective structure  $\mu$  and 
a packing  $P$. 
  
The problem then is to find suitable conditions 
which must be satisfied by the cross ratio parameters. 
To do this, 
for each map  ${\bf c} : E(\tau) \to {\mathbb R}_+$  which 
does not necessarily come from a circle packing, 
we assign to each edge  $e$  an associated matrix 
$A={\left(\begin{array}{cc} 0 & 1 \\ -1 & {\bf c}(e) \\ \end{array}\right)} 
\in SL_2({\mathbb R})$.  
Furthermore, 
for each vertex $v \in V(\tau)$  with valence $n$, 
one reads off the edges incident to $v$ in 
the clockwise direction and form a word  $W_v$  of length  $n$  from 
the corresponding associated matrices   
($W_v$ is defined up to cyclic permutation).  
The main general result is the following: 

\smallskip

\begin{lemma*}[Lemma \ref{Lem:MainLemma}] 
Suppose  ${\bf c}_P$  is the cross ratio parameter of 
the packing  $P$  on  $(\Sigma_g, \mu)$  
with nerve  $\tau$    
and let  $W_v$  be the associated word of 
the associated matrices for each  $v \in V(\tau)$.  
Then, we have  
\begin{enumerate}
\item[$(1)$] 
	$W_v=-I$,  where $I$ is the identity matrix, and 
\item[$(2)$] 
	if the length of  $W_v$  is  $n$,  
	then every subword of $W_v$ of length $\leq n-1$  is admissible, 
	and every subword of length $\leq n-2$ are strictly admissible,  
	(the admissibility and strict admissibility are defined by some 
	inequalities in the cross ratio type invariants which will be 
	precisely defined in \S 2). 
\end{enumerate}
Conversely, 
if a map  ${\bf c} : E(\tau) \to {\mathbb R}_+$  satisfies  $(1)$  and  $(2)$, 
then there  is a unique projective structure  $\mu$  on  
$\Sigma_g$  together 
with a circle packing  $P$  with nerve isotopic to  $\tau$  
such that the cross ratio parameter  ${\bf c}_P$  is equal to  ${\bf c}$. 
\end{lemma*} 

\medskip

We call the set of all maps in ${\mathbb R}^{E(\tau)}$ satisfying 
conditions (1) and (2) for each vertex  $v \in V(\tau)$  in 
the main lemma the cross ratio parameter space 
and denote it by
\begin{equation*} 
	{\mathcal C}_{\tau} = 
	\{ {\bf c} : E(\tau) \to {\mathbb R}_+ \, \vert \, 
	{\bf c} \; \; \text{satisfies (1) and (2) 
		for each vertex} \}. 
\end{equation*}  
By this lemma, 
${\mathcal C}_{\tau}$  can be identified with 
the set of all pairs  $(\mu, P)$  and hence 
contains a special element obtained by the Andreev-Thurston theorem. 
  
Clearly, ${\mathcal C}_{\tau}$ is a real semi-algebraic set,
however this does not say much about its true structure.  
It is a priori possible that ${\mathcal C}_{\tau}$ has singularities 
or may even reduce to a single point 
represented by the Andreev-Thurston solution. 
Also, each point in ${\mathcal C}_{\tau}$ determines not just 
a projective structure but a circle packing on the surface 
with nerve $\tau$. 
It is quite possible that there is 
a non-trivial family of different circle packings with 
nerve $\tau$ on the surface with a fixed projective structure
(that is, the circle packings are not rigid).  

Such problems can be well formulated using the map 
\begin{equation*} 
	f : {\mathcal C}_{\tau} \to {\mathcal P}_g 
\end{equation*} 
defined by forgetting the packing  $P$. 
For example, 
the rigidity question is equivalent to the  injectivity of  $f$.  

In \S 3, 
we will count the formal dimension of  ${\mathcal C}_{\tau}$  which 
is just the number of variables minus the number of equations 
which appear in the definition, and show that it is  $6g-6$  when  $g \geq 1$.  
Then we show using the deformation theory of Kleinian groups 
that this turns out to give the correct dimension at least near 
the Andreev-Thurston solution when   $g \geq 2$.   
When  $g = 1$, 
we will see by  hyperbolic Dehn surgery theory that 
the formal dimension count overlooks some dependence in the equations. 
The results in \S 3 can be summarized by 

\smallskip 

\begin{lemma*}[Lemma \ref{Lem:NbdAT2}, Lemma \ref{Lem:NbdAT1}] 
The Andreev-Thurston solution 
has a neighborhood in  ${\mathcal C}_{\tau}$  
which is homeomorphic to the euclidean space 
of dimension  $6g-6$  when  $g \geq 2$  and  $2$  when  $g = 1$, 
and to which the restriction of the forgetting map  $f$  is 
an embedding.  
\end{lemma*} 

\medskip

We now state our main result which describes 
a global picture of  ${\mathcal C}_{\tau}$  for 
packings by one circle. 
In this case, 
$\tau$  has exactly one vertex.  
Note that a graph with one vertex which triangulates  $\Sigma_g$  is 
necessarily covered by a simple graph in the universal cover. 

\smallskip

\begin{maintheorem}[Lemma \ref{Lem:Cell}, 
Lemma \ref{Lem:Rigid}\label{Thm:MainTheorem}]  
Suppose that the triangulation  $\tau$  of  $\Sigma_g$  ($g \geq 2$)  
has exactly one vertex. 
Then 
\begin{itemize}
\item[$(i)$] 
	${\mathcal C}_{\tau}$ is homeomorphic to 
	the euclidean space of dimension  $6g-6$, and 
\item[$(ii)$] 
	the forgetting map  $f : {\mathcal C}_{\tau} \to {\mathcal P}_g$  
	is an embedding. 
\end{itemize}
\end{maintheorem} 

\medskip
\noindent {\bf Remarks:}
1. When  $g = 1$, 
the same conclusion  (replacing the dimension by  $2$)
can be deduced from a previous work by Mizushima  \cite{Miz}. 
His methods, different from ours,  
provides much more  information about the deformation space, as we see in the third remark, 
but do not generalize easily to the higher genus case. 
Our method works also for this case, and it is  discussed in the appendix.
\smallskip

2. We conjecture that the main theorem holds for any triangulation $\tau$ of $\Sigma_g$  which is covered by a simple graph in the universal cover.  
\smallskip 

3. In \cite{Miz}, 
Mizushima showed that the composition  $u \circ f$  of the 
forgetting map and the uniformization is a homeomorphism when $g=1$ 
(stated in a slightly different way).  
We conjecture that this holds for genus $g \geq 2$ as well.  
If true, this would give a uniformization theorem 
that every complex structure can be uniformized by 
a projective structure circle-packed by one circle.  
However, again, the methods of \cite{Miz} do not generalize.  
In the torus case, it is relatively easy to uniformize complex affine 
structures on the torus to Euclidean structures via the $\log$ function.  
The general uniformization map is much more complicated and 
so different methods would have to be employed.

\medskip

We here outline the proof of the main theorem 
and show how the rather strong condition 
that  $\tau$ has exactly one vertex plays a crucial role.  
In this case, $\tau$  has $6g-3$ edges 
so that there are $6g-3$ variables. 
Condition (1) of the main lemma then corresponds to 3 relations 
on these variables,   
To show non-singularity of  ${\mathcal C}_{\tau}$,  
we will find a specific set of three variables such that the 
remaining  $6g-6$  variables are free but the chosen three 
are uniquely determined by these free variables.  
Then we will see that the range of free variables 
consists of a convex set.   
This is enough to see that  ${\mathcal C}_{\tau}$  is 
homeomorphic to the euclidean space of dimension  $6g-6$. 

To see the injectivity of the forgetting map, 
suppose that  ${\bf c}$ and ${\bf c'}$ are two cross ratio parameters  which correspond to the same projective structure on $\Sigma_g$ (but possibly different circle packings).  Then we would like to show that ${\bf c}={\bf c'}$ which implies rigidity of the circle packing. A priori, we only know that the holonomy representations corresponding to ${\bf c}$ and ${\bf c'}$ are conjugate since they give rise to the same projective structure. We show that  the holonomy representation of the specific set of three side-pairings which appear in showing nonsingularity of  ${\mathcal C}_{\tau}$  must be actually the same for ${\bf c}$ and ${\bf c'}$, and not just conjugate. We then show that for such a triple, the holonomy image is always non-elementary so that the holonomy representation of the full fundamental group $\pi_1(\Sigma_g)$ is exactly the same for both ${\bf c}$ and ${\bf c'}$ and hence ${\bf c}={\bf c'}$.

\medskip 
\noindent {\bf Remark:} 
The argument above for the injectivity of the forgetting 
map  $f$  implies a rather remarkable conclusion that 
the composition of  $f$  with the map  
$\chi : {\mathcal P}_g \to Hom(\pi_1(\Sigma_g), 
	PSL_2({\mathbb C}))/ PSL_2({\mathbb C})$  
which assigns to each projective structure the conjugacy class of 
its holonomy representation is injective.  
This contrasts with the fact that  $\chi$, which  
is a local homeomorphism, is in fact never injective, 
see for instance  \cite{Hej}.

\medskip

The rest of this paper is organized as follows. In \S 2, we give some preliminaries, introduce the cross ratio and associated matrix, and prove the general result that the cross ratio parameter space associated to $\tau$ is a real semi-algebraic set. We also show how to obtain the holonomy representation from the cross ratios. In \S 3, we prove Lemma 2 using well developed machinery from  3-dimensional deformation theory.  In \S 4, we study  the cross ratio parameter space ${\mathcal C}_{\tau}$ when $\tau$ has one vertex, which corresponds to a circle packing by one circle. In particular, we show that  ${\mathcal C}_{\tau}$ can be parameterized by coordinates lying in a full convex subset of ${\mathbb R}^{6g-6}$. In \S 5, we prove  that the circle packing is rigid for the one circle packing, hence completing the proof of our main theorem. Finally, we work out in detail the case of the one circle packing on the torus ($g=1$) in the appendix, where although the results are similar, the methods are slightly different.

%
%

\section{Cross Ratio Parameter} 


\subsection{Preliminaries}

A projective structure on a closed orientable surface $\Sigma_g$ of genus  $g$   is a maximal collection of local charts modeled on the Riemann sphere  $\hat{\mathbb C}$ such that locally, the transition functions are restrictions of $\gamma \in PSL_2({\mathbb C})$. Another point of view is that a projective structure is a structure modeled on the boundary at $\infty$ of the three-dimensional hyperbolic space ${\mathbb H}^3$ under the action of the isometries of  ${\mathbb H}^3$. Associated to such a projective structure is a developing map $dev: \widetilde {\Sigma}_g \longrightarrow \hat{\mathbb C}$, defined up to composition with elements of $PSL_2({\mathbb C})$ and a holonomy representation $\rho: \pi_1(\Sigma_g) \longrightarrow PSL_2({\mathbb C})$ defined up to conjugation by elements of $PSL_2({\mathbb C})$ such that the holonomy representation is equivariant with respect to the developing map.  

If $\mu$ is a projective structure on $\Sigma_g$,  a disk in $(\Sigma_g, \mu)$ is the closure of a simply connected region ${\mathcal D}^i$ in $\Sigma_g$ such that the lifts of the closure  ${\mathcal D}$ in $\widetilde{\Sigma}_g$ are  mapped homeomorphically to closed disks in $\hat{\mathbb C}$ by the developing map. This makes sense since circles are preserved by the linear fractional transformations $\gamma \in PSL_2({\mathbb C})$. 

\smallskip

\begin{definition} 
A circle packing on $(\Sigma_g, \mu)$ is a collection of closed disks $\{{\mathcal D_i}\}$ in $\Sigma_g$ such that the interiors of any two distinct disks are disjoint and $\Sigma_g \backslash \cup_i {\mathcal D_i}$ consists of a finite number of triangular interstices each bounded by three circular arcs.
\end{definition} 
\smallskip

\begin{definition} 
For a circle packing on $\Sigma_g$, we assign a vertex to each circle and  an edge joining two vertices for each tangency point between two circles. The  graph $\tau$ on $\Sigma_g$ obtained thus triangulates $\Sigma_g$  and this is called the nerve of the circle packing. 
\end{definition} 

\smallskip


\subsection{Cross ratio and associated matrix}

We  define a cross ratio type invariant on the edges of $\tau$, which has the property that if an edge separating two interstices is assigned a given positive real number, then knowing the developing image of one of the interstices determines the developing image of the other interstice. We start with a collection of 4 circles in $\hat{\mathbb C}$ whose nerve is a graph consisting of 4 vertices and 5 edges such that the edges bound two triangles with a common edge. The cross ratio type invariant (henceforth the cross ratio for short) is defined on the common edge (see Figure \ref{Fig:Configuration}a).

\smallskip

\begin{proposition}\label{Prop:CrossRatio}
Let $C_1,C_2,C_3,C_4$ be 4 circles in $\hat{\mathbb C}$ such that $C_1$ is tangent to $C_2$, $C_3$ and $C_4$, $C_2$ is tangent to $C_3$, and $C_3$ is tangent to $C_4$ (see Figure \ref{Fig:Configuration}b). Denote the tangency points between $C_i$ and $C_j$ by $p_{ij}$. We also further assume that $p_{12}$, $p_{13}$ and $p_{14}$ lie in a clockwise direction on $C_1$. Then the cross ratio $(p_{14}, p_{23},p_{12},p_{13})=x \sqrt{-1}$, where $x \in {\mathbb R}_+$.
\end{proposition}
 \smallskip

\noindent {\it Proof:} We can find a linear fractional transformation $\gamma \in PSL_2({\mathbb C})$ mapping $p_{23}$ to $1$, 
$p_{12}$ to 0 and $p_{13}$
to $\infty$ so that $\gamma(p_{14})$$=(p_{14}$, $p_{23}$, $p_{12}$, $p_{13})$ using the  definition of the cross ratio of 4 points given in Ahlfors \cite{Ahl}. Clearly, $\gamma(C_1)=\{ z \in {\mathbb C}| Re(z)=0\} \cup \{\infty\}$, $\gamma(C_2)=\{ z \in {\mathbb C}|~|z-1/2|=1/2\}$ and $\gamma(C_3)=\{z \in {\mathbb C}|Re(z)=1\}\cup\{\infty\}$. It follows that $\gamma(p_{14})$ is a purely imaginary number since it lies on $\gamma(C_1)$. The condition $p_{12}$, $p_{13}$ and $p_{14}$ lie in a clockwise direction on $C_1$ now implies that $\gamma(p_{14})=x\sqrt{-1}$ where $x>0$, see Figure \ref{Fig:CrossRatio}. 
\qed

\begin{figure}[ht]
  \begin{center}
    \begin{picture}(300,140)
        \put(0,0){\scalebox{0.7}{\includegraphics{Configuration-a.eps}}}
        \put(160,0){\scalebox{0.7}{\includegraphics{Configuration-b.eps}}}
\put(64,4){(a)}
\put(191,81){$p_{23}$}
\put(191,59){$p_{12}$}
\put(228,78){$p_{13}$}
\put(257,81){$p_{34}$}
\put(257,58){$p_{14}$}
\put(243,118){$C_3$}
\put(293,50){$C_4$}
\put(246,20){$C_1$}
\put(155,84){$C_2$}
\put(224,4){(b)}
    \end{picture}
  \end{center}
  \caption{}
\label{Fig:Configuration}
\end{figure}

\begin{figure}[ht]
  \begin{center}
    \begin{picture}(140,160)
      \put(0,0){\scalebox{0.7}{\includegraphics{CrossRatio.eps}}}
\put(-21,109){$ix=\gamma(p_{14})$}
\put(29,148){$\gamma(C_1)$}
\put(72,148){$\gamma(C_3)$}
\put(85,109){$\gamma(p_{34})$}
\put(127,143){$\gamma(p_{13})=\infty$}
\put(85,20){$\gamma(p_{23})=1$}
\put(47,56){$\gamma(C_2)$}
\put(-15,20){$\gamma(p_{12})=0$}
    \end{picture}
  \end{center}
  \caption[]{}
\label{Fig:CrossRatio}
\end{figure}

\begin{definition}[Cross Ratio of Configuration]\label{Def:CrossRatio}
Given a configuration of 4 circles as defined in proposition \ref{Prop:CrossRatio} above, the cross ratio of the configuration, or alternatively, of the interior edge  of the nerve of the configuration  is the positive real number $x$ such that  $(p_{14}, p_{23},p_{12},p_{13})=x \sqrt{-1}$. 
\end{definition} 

\smallskip

\noindent {\bf Remarks:}
1. It is easy to see that the cross ratio of $(C_1,C_2,C_3,C_4)$ is the same as that for $(C_3,C_4,C_1,C_2)$ so that the cross ratio is defined on the undirected interior edge.
\smallskip 

2. The cross ratio  is clearly an invariant of the configuration under the action of $PSL_2({\mathbb C})$. It measures the distance between $C_2$ and $C_4$ in the following sense: Consider the circles $C_2$ and $C_4$ as the boundaries at infinity of two hyperbolic planes in ${\mathbb H}^3$.  If the two circles are disjoint or tangent, the planes are at some non-negative distance $d$ from each other. Then the cross ratio  $x=\cosh(d/2)$. If the two circles intersect at angle $\theta$, then $x=\cos (\theta/2)$. It follows that $C_2$ and $C_4$ are disjoint if and only if $x>1$.
\smallskip

3. The cross ratio, together with $C_1$, $C_2$ and $C_3$ determines $C_4$. Alternatively, we can say that knowing the position of one interstice and the cross ratio determines the neighboring one.  
\medskip

\begin{definition}[Associated Matrix]\label{Def:AssoMatrix}
To each configuration of 4 circles having the cross ratio  $x$  
is associated the matrix 
$A={\left(\begin{array}{cc} 0 & 1 \\ -1 & x \\ \end{array}\right)}$.  
\end{definition} 

We note that $A\in SL_2({\mathbb R})$ rather than $PSL_2({\mathbb R})$. 
This is important as it will allow us to keep track of overlapping in a configuration of one circle surrounded by several circles.

We next explain the role of the associated matrix $A$. Basically, if we put the configuration of 4 circles into a standard position by some linear fractional transformation, then $A$ maps $C_1$ to itself and one interstice to the other. As much of the subsequent arguments depend on putting a configuration into some standard position, we start with the following definition.

\smallskip

\begin{definition}[Standard Interstice]\label{Def:Interstice}
The standard interstice is defined to be the interstice with vertices at $0$, $\sqrt{-1}$ and $\infty$. It is bounded by the 3 circles (or straight lines) $C_1=\{z\in \hat{\mathbb C}:Im(z)=0\}$, $C_2=
\{z\in \hat{\mathbb C}:Im(z)=1\}$ and $C_3=\{z\in \hat{\mathbb C}: ~|z-\sqrt{-1}/2|=1/2\}$. We denote the standard interstice by $\mathcal{I}_s$.
\end{definition} 

\smallskip

\begin{proposition}\label{Prop:MShift}
Let $C_1$, $C_2$, $C_3$ and $C_4$ be a configuration of 4 circles as in proposition \ref{Prop:CrossRatio} with cross ratio $x$ and associated matrix 
$A={\left(\begin{array}{cc} 0 & 1 \\ -1 & x \\ \end{array}\right)}$ 
as defined in definition \ref{Def:AssoMatrix}. Furthermore, suppose that $C_1$, $C_2$ and $C_3$ are as defined in definition \ref{Def:Interstice} so they bound the standard interstice (see Figure \ref{Fig:Interstice}). Then $A(C_1)=C_1$, $A(C_2)=C_3$, $A(C_3)=C_4$. In other words, $A$ maps the interstice $\mathcal{I}_s$ to the other interstice of the configuration, with the boundary adjacent to $C_1$ being mapped to the boundary adjacent to $C_1$.
\end{proposition} 

\smallskip

\noindent {\it Proof:}
First note that by a simple computation, since the cross ratio is $x$, $p_{14}=1/x$. Clearly, $A(C_1)=C_1$ since entries of $A$ are real. Next  $A(C_2)=C_3$ since $A(\infty)=0$ and $A(x+\sqrt{-1})=\sqrt{-1}$ and $A$ maps tangency points to tangency points. Finally, $A(C_3)=C_4$ since $A(0)=1/x$. \qed

\begin{figure}[ht]
  \begin{center}
    \begin{picture}(140,50)
\put(0,-20){
      \put(0,0){\scalebox{0.7}{\includegraphics{Interstice.eps}}}
\put(146,15){$C_1$}
\put(146,48){$C_2$}
\put(33,33){$C_3$}
\put(57,27){$C_4$}
\put(4,27){$\mathcal{I}_s$}
}
    \end{picture}
  \end{center}
  \caption[]{}
\label{Fig:Interstice}
\end{figure}

For convenience, we state in the next proposition various configurations of four circles $C_1,C_2,C_3,C_4$ all of which have cross ratio  $x$.

\smallskip

\begin{proposition}\label{Prop:SomeView}
The  configurations of four circles $C_1,C_2,C_3,C_4$ defined in each of the following cases all have cross ratio  $x$ (see Figure \ref{Fig:SomeView}): 
\begin{itemize} 
\item[(i)] 
	$C_1=\{z\in \hat{\mathbb C}:Im(z)=0\}$, $C_2=\{z\in \hat{\mathbb C}: ~|z-(x+\sqrt{-1}/2)|=1/2\}$,  $C_3=\{z\in \hat{\mathbb C}:Im(z)=1\}$, $C_4=\{z\in \hat{\mathbb C}: ~|z-\sqrt{-1}/2|=1/2\}$.
\item[(ii)] 
	$C_1=\{z\in \hat{\mathbb C}:Im(z)=0\}$, $C_2=\{z\in \hat{\mathbb C}:Im(z)=1\}$, $C_3=\{z\in \hat{\mathbb C}: ~|z-\sqrt{-1}/2|=1/2\}$, $p_{14}=1/x$.
\item[(iii)] 
	$C_1=\{z\in \hat{\mathbb C}:Im(z)=0\}$, $p_{12}<p_{13}<p_{14}$, ${\displaystyle \sqrt{\frac{r_3}{r_2}}+\sqrt{\frac{r_3}{r_4}}=x}$, where $r_i$ is the radius of $C_i$, $i=2,3,4$.
\end{itemize} 
\end{proposition} 

\medskip

\noindent {\it Proof:} The proof follows from some elementary calculations and will be omitted.
\qed 

\begin{figure}[ht]
  \begin{center}
    \begin{picture}(320,200)
      \put(0,110){\scalebox{0.7}{\includegraphics{SomeView-1.eps}}}
      \put(180,110){\scalebox{0.7}{\includegraphics{SomeView-2.eps}}}
      \put(0,15){\scalebox{0.7}{\includegraphics{SomeView-3.eps}}}
\put(60,90){(i)}
\put(-17,123){$C_1$}
\put(106,143){$C_2$}
\put(-17,159){$C_3$}
\put(34,144){$C_4$}
\put(32,114){$0$}
\put(23,170){$\sqrt{-1}$}
\put(103,114){$x$}

\put(240,90){(ii)}
\put(165,123){$C_1$}
\put(165,159){$C_2$}
\put(216,146){$C_3$}
\put(241,138){$C_4$}
\put(212,114){$0$}
\put(203,170){$\sqrt{-1}$}
\put(240,114){$p_{14}=\frac{1}{x}$}

\put(0,-5){
\put(57,0){(iii)}
\put(-14,31){$C_1$}
\put(33,28){$p_{12}$}
\put(63,28){$p_{13}$}
\put(94,28){$p_{14}$}
\put(39,47){$r_2$}
\put(68,45){$r_3$}
\put(100,47){$r_4$}
\put(39,75){$C_2$}
\put(66,68){$C_3$}
\put(100,75){$C_4$}
}
    \end{picture}
  \end{center}
  \caption[]{}
\label{Fig:SomeView}
\end{figure}

\medskip


\subsection{Admissibility}
 
In the next set of results, we  consider a configuration of circles with a central circle $\overline{C}$ and $n+1$ circles $C_0,\cdots,C_n$ surrounding it such that  
\begin{itemize}  
\item[(i)] 
	$\overline{C}$ is tangent to $C_i$ for $i=0,\cdots,n$; 
\item[(ii)]
	$C_i$ is tangent to $C_{i-1}$ and $C_{i+1}$, for $i=1, \cdots,n-1$; 
\item[(iii)]
	The points of tangency $p_i$ between $\overline{C}$ and $C_i$ are arranged in clockwise direction around $\overline{C}$ (see Figure \ref{Fig:Admissible}a).
\end{itemize}

Let $\tau$ be the nerve of this configuration, with vertices $\overline{v}, v_0, \cdots, v_n $ corresponding to the circles $\overline{C}, C_0,\cdots,C_n$ respectively.  Denote by $e_i$ the edge joining $\overline{v}$ to $v_i$ and $e_{\{i,i+1\}}$ the edge  joining $v_i$ to $v_{i+1}$.

\medskip 

\begin{definition}[Cross Ratio Vector and Associated Word]\label{Def:Vector}
For such a configuration,  let $x_i$ be the cross ratio of the edge $e_i$  and $A_i$ be the associated matrix, for $i=1, \cdots, n-1$. ${\bf x}=(x_1,x_2, \cdots ,x_{n-1})$ will be called the {\it cross ratio vector} of the configuration and $W_{\bf x}=A_1A_2 \cdots A_n$ the {\it word of associated matrices} of the configuration. 
\end{definition} 

We have the following result:

\medskip

\begin{proposition}\label{Prop:SurroundConfig}
Suppose that $\{\overline{C}, C_0,\cdots,C_n \}$ is a configuration of circles as defined above. Let $I_i$ be the interstice bounded by the circles $\overline{C}$, $C_i$ and $C_{i+1}$. We say that the configuration is in standard position if $I_0$ is the standard interstice $\mathcal{I}_s$ and $\overline {C}$ is the real line (see Figure \ref{Fig:Admissible}b). Let $W_k=A_1A_2\cdots A_k$ \, for $k=1, \cdots n-1$. If the configuration is in standard position, then $W_k(I_0)=I_k$ with the side of $I_0$ bounded by $\overline{C}$ being mapped to the side of $I_k$ bounded by $\overline{C}$, so that $W_k(p_0)=p_{k}$, $W_k(p_1)=p_{k+1}$.
\end{proposition} 

\medskip

\noindent {\it Proof: } We prove by induction on $k$. For $k=1$, the result is just proposition \ref{Prop:MShift}. Suppose the result is true for $k$. Then $ W_k A_{k+1} W_k^{-1}(I_k)=I_{k+1} \Longrightarrow W_k A_{k+1} W_k^{-1}W_k(I_0)=I_{k+1} \Longrightarrow W_{k+1}(I_0)=I_{k+1}$. Also, since $W_k \in SL_2({\mathbb R})$, $W_k$ maps $\overline{C}$ to itself so the second condition holds. \qed

\begin{figure}[ht]
  \begin{center}
    \begin{picture}(320,140)
      \put(0,0){\scalebox{0.7}{\includegraphics{Admissible-a.eps}}}
      \put(160,10){\scalebox{1.0}{\includegraphics{Admissible-b.eps}}}
\put(60,-10){(a)}
\put(66,65){$\overline{C}$}
\put(67,101){$p_0$}
\put(94,87){$p_1$}
\put(99,65){$p_2$}
\put(34,64){$p_n$}
\put(85,128){$C_0$}
\put(121,97){$C_1$}
\put(120,62){$C_2$}
\put(12,51){$C_n$}

\put(260,-10){(b)}
\put(0,10){
\put(185,7){$p_1$}
\put(212,7){$p_2$}
\put(230,7){$p_3$}
\put(288,7){$p_n$}
\put(186,28){$C_1$}
\put(296,27){$C_n$}
\put(147,11){$\overline{C}$}
\put(145,46){$C_0$}
\put(135,31){$I_0=\mathcal{I}_s$}
}
    \end{picture}
  \end{center}
  \caption[]{}
\label{Fig:Admissible}
\end{figure}

In the configuration of one central circle and $n+1$ surrounding circles defined in proposition \ref{Prop:SurroundConfig}, it is possible that some of the surrounding circles $C_0,\cdots, C_n$ not adjacent to each other may overlap. This is permissible from the definition, in general, from the point of view of circle packings on surfaces with projective structures,  this would imply that the developing map from $\widetilde{\Sigma}_g$ to $\hat{\mathbb C}$ is not a homeomorphism onto its image. More importantly, it may be possible that boundaries of the interstices lying on the central circle ${\overline{C}}$ overlap as they wrap around  ${\overline{C}}$. We are interested in the case where they do not overlap, that is, they do not go more than once around ${\overline{C}}$. In standard position, this just means that the points of tangency $p_i$ between $\overline{C}$ and $C_i$ satisfy the inequalities
$$p_0=-\infty<p_1=0 <p_2 <\cdots <p_n.$$

\medskip

\begin{definition}[Admissibility]\label{Def:Admissibility}
The circles $\{\overline{C}, C_0,\cdots,C_n\}$ as given in proposition \ref{Prop:SurroundConfig} is   {\it admissible} if in standard position, the points of tangency  $p_i$ between $\overline{C}$ and $C_i$ satisfy $$p_0=-\infty<p_1=0 <p_2 <\cdots <p_n.$$  It is  {\it strictly admissible} if $p_n \neq \infty$, that is, if in addition, $p_n \neq p_0$.  The $(n-1)$-tuple  ${\bf x}=(x_1,\cdots,x_{n-1}) \in {\mathbb R}_+^{n-1}$ is said to be  (strictly) admissible if it is the cross ratio vector of a (strictly) admissible configuration. Similarly, a word $W_{\bf x}=A_1A_2 \cdots A_{n-1}$ is  (strictly) admissible if it is the word of associated matrices of a (strictly) admissible configuration.
\end{definition} 

\medskip

 Geometrically, a strictly admissible configuration is one where the boundaries of the interstices lying on the central circle do not overlap, an admissible configuration is one where there is either no overlap or overlap at only one point. It is clear from the definition that a proper subword of an admissible word is always strictly admissible.

\medskip
The condition of admissibility of $(x_1,\cdots,x_{n-1}) \in {\mathbb R}_+^{n-1}$ can be easily translated into a condition on the subwords of the word $W_{\bf x}=A_1 A_2 \cdots A_{n-1}$ in the associated matrices as seen below:

 \smallskip

\begin{proposition}\label{Prop:CharacterizeAdm}
Let $(x_1,\cdots,x_{n-1}) \in {\mathbb R}_+^{n-1}$ and 
$A_i={\left(\begin{array}{cc} 0 & 1 \\ -1 & x_i\\ \end{array}\right)}$ $\in SL_2({\mathbb R})$,  $i=1,\cdots,n-1$.
Then $(x_1,\cdots,x_{n-1}) \in {\mathbb R}_+^{n-1}$ is admissible if and only if for all subwords 
$W^i_k={\left(\begin{array}{cc} a_{ik} & b_{ik} \\ c_{ik} & d_{ik} \\ \end{array}\right)}  =A_iA_{i+1} \cdots A_k$ of $W_{\bf x}$,  $a_{ik} \leq 0$ (with equality only when $i=k$), $b_{ik} > 0$, $c_{ik} < 0$, and $d_{ik} \geq 0$. Furthermore, the last inequality is always strict for all proper subwords, and if it is also strict for the word $W_{\bf x}$, then the word is strictly admissible. 
\end{proposition}

\medskip

\noindent {\it Proof:} For convenience, we consider only the case of strict admissibility. The remaining case is easily deduced.

$(\Longleftarrow)$
Consider the configuration $\{\overline{C}, C_0,\cdots,C_n\}$ in standard position with cross ratio vector $(x_1,\cdots,x_{n-1}) \in {\mathbb R}_+^{n-1}$. Recall that $p_i$ is the point of tangency of $\overline{C}$ with $C_i$ and $p_0=-\infty$, $p_1=0$ since the configuration is in standard position.   To simplify notation, let $W^1_k={\left(\begin{array}{cc} a_k & b_k \\ c_k & d_k \\ \end{array}\right)}$. 
By proposition \ref{Prop:SurroundConfig}, $W^1_k(0)=b_k/d_k=p_{k+1}$, $W^1_k(-\infty)=a_k/c_k=p_k$.
Since $b_k,d_k > 0$, $p_{k+1}>0$ for $k=1, \cdots, n-1$. Furthermore, since $a_kd_k-b_kc_k=1$, $p_k=p_{k+1}+1/(c_kd_k)<p_{k+1}$ since $c_k<0,d_k>0$.  
Hence,  
$$p_0=-\infty <p_1=0<p_2<p_3 \cdots <p_n < \infty $$ so that the configuration is strictly admissible. 

\medskip
$(\Longrightarrow)$
We first note that if $\{\overline{C}, C_0, C_1, \cdots, C_n\}$ is strictly admissible, then so is the sub-configuration $\{\overline{C}, C_i, C_{i+1}, \cdots, C_n\}$, so $(x_1,\cdots ,x_{n-1})$ is admissible implies $(x_i, x_{i+1}, \cdots,x_{n-1})$ is admissible for $0<i \leq n-1$. Hence it suffices to prove that $W^1_k={\left(\begin{array}{cc} a_k & b_k \\ c_k & d_k \\ \end{array}\right)}$ satisfies the conditions stated in the proposition for $k=1, \cdots, n-1$. We do this by induction on $k$. For $k=1$, $W^1_1={\left(\begin{array}{cc} 0 & 1 \\ -1 & x_1\\ 
\end{array}\right)}$ which clearly satisfies the conditions. Now suppose that $W^1_k$ satisfies the conditions. Then $$W^1_{k+1}={\left(\begin{array}{cc} a_{k+1} & b_{k+1} \\ c_{k+1} & d_{k+1}\\ 
\end{array}\right)}=
{\left(\begin{array}{cc} a_k & b_k \\ c_k & d_k\\ 
\end{array}\right)}{\left(\begin{array}{cc} 0 & 1 \\ -1 & x_{k+1}\\ 
\end{array}\right)}$$
$$={\left(\begin{array}{cc} -b_k & a_k+b_kx_{k+1} \\ -d_k & c_k+d_kx_{k+1}\\ 
\end{array}\right)}.$$
Hence $a_{k+1}=-b_k<0$, $c_{k+1}=-d_k<0$. Furthermore, since the configuration is admissible, ${\displaystyle \frac{b_{k+1}}{d_{k+1}}=p_{k+1}>p_k=\frac{a_{k+1}}{c_{k+1}}>0}$.
Using the same argument as in the previous case, we have $$a_{k+1}d_{k+1}-b_{k+1}c_{k+1}=1$$
$$ \Longrightarrow \frac{a_{k+1}}{c_{k+1}}-\frac{b_{k+1}}{d_{k+1}}=\frac{1}{c_{k+1}d_{k+1}}<0. $$
Since $c_{k+1}<0$, this implies $d_{k+1}>0$ and hence $b_{k+1}>0$. 
\qed

\medskip
\begin{example} 
{\em 
It is easily checked that the vector $(\sqrt{2},\sqrt2,\sqrt2,\sqrt2)$ and corresponding word ${\tiny {\left(\begin{array}{cc} 0 & 1 \\ -1 & \sqrt2\\ \end{array}\right)}{\left(\begin{array}{cc} 0 & 1 \\ -1 & \sqrt2\\ \end{array}\right)}{\left(\begin{array}{cc} 0 & 1 \\ -1 & \sqrt2\\ \end{array}\right)}{\left(\begin{array}{cc} 0 & 1 \\ -1 & \sqrt2\\ \end{array}\right)}}$ are admissible (but not strictly admissible) whereas the vector \\
$(\sqrt{2},\sqrt2,\sqrt2,\sqrt2,\sqrt2)$ and corresponding word are not admissible.
}
\end{example}
\medskip

The next proposition tells us when a strictly admissible word $W$ remains admissible when we  multiply it by  an associated matrix either on the left or right, or both.

\medskip

\begin{proposition}\label{Prop:ByX}
Suppose $W={\left(\begin{array}{cc} a & b \\ c & d\\ \end{array}\right)}=A_1A_2 \cdots A_n$ is a strictly admissible word and $X={\left(\begin{array}{cc} 0 & 1 \\ -1 & x\\ 
\end{array}\right)}$. Then 
\begin{itemize}
\item[(i)] 
	$XW$ is admissible if and only if $x \geq b/d$ ; 
\item[(ii)] 
	$WX$ is admissible if and only if $x \geq -c/d$ ; 
\item[(iii)] 
	$XWX$ is admissible if  and only if $${\displaystyle x \geq \frac{b-c+\sqrt{(b+c)^2+4}}{2d} }$$
with strict admissibility if the inequalities are strict.
\end{itemize}
\end{proposition} 

\medskip

\noindent {\it Proof:} (i) Let $\{\overline{C}, C_0, \cdots,C_{n+1}\}$ be the configuration in standard position corresponding to $W$. The configuration corresponding to $AW$ is obtained by adding another circle $C_{-1}$ tangent to $\overline{C}$ and $C_0$. This configuration is admissible if and only if $p_{-1} \geq p_{n+1}$, see Figure \ref{Fig:byA}, with strict admissibility if the inequality is strict. By proposition \ref{Prop:SomeView}, $p_{-1}=x$, by proposition \ref{Prop:SurroundConfig}, $p_{n+1}=b/d$ so (i) follows. 

(ii) The proof is similar to that of proposition \ref{Prop:CharacterizeAdm} and will be omitted. 

(iii) Let
$\{\overline{C}, C_0, \cdots,C_{n+1}\}$ be the configuration in standard position corresponding to $W$ and let $\{\overline{C}, C_{-1}, C_0, \cdots,C_{n+1}, C_{n+2}\}$ be the augmented configuration corresponding to $XWX$. Clearly, the configuration is admissible if and only if $p_{-1} \geq p_{n+2}$ and conditions (i) and (ii) are satisfied. Since $p_{-1}=x$ and ${\displaystyle p_{n+2}=\frac{a+bx}{c+dx}}$. This is easily shown to be equivalent to $x \geq \alpha$ where $\alpha$ is the larger of the two roots of the equation $dx^2+(c-b)x-a=0$ which is equivalent to (iii) using $ad-bc=1$.
\qed

\begin{figure}[ht]
  \begin{center}
    \begin{picture}(200,60)
      \put(0,-10){\scalebox{1}{\includegraphics{byA.eps}}}
\put(0,-10){
\put(28,7){$p_1$}
\put(53,7){$p_2$}
\put(70,7){$p_3$}
\put(125,7){$p_{n+1}$}
\put(175,7){$p_{-1}$}
\put(28,28){$C_1$}
\put(54,23){$C_2$}
\put(130,31){$C_{n+1}$}
\put(170,28){$C_{-1}$}
\put(-11,13){$\overline{C}$}
\put(-15,46){$C_0$}
}
    \end{picture}
  \end{center}
  \caption[]{}
\label{Fig:byA}
\end{figure}
%
%


\subsection{Convexity} 

The next result shows the convexity of the admissibility condition, that is, fixing some positive $n$, the set of all strictly admissible cross ratio vectors of length $n$ is convex.

Let $C_{SA}^n\subset {\mathbb R}^n$ be the set of all strictly admissible cross ratio vectors of length $n$ and let ${\bf x}=(x_1,\cdots,x_n) \in C_{SA}^n$.   Suppose $W_{\bf x}=A_1A_2\cdots A_n={\left(\begin{array}{cc} a & b \\ c & d\\ \end{array}\right)}$ is the word of associated matrices and let $p=b/d$. 
We now regard the entries  $a$, $b$, $c$, $d$ of $W_{\bf x}$ and the ratio  $p$ as  real valued functions of ${\bf x}$. 

\medskip

\begin{proposition}\label{Prop:Convex}
We have the following: 
\begin{itemize}
\item[(i)] 
	$C_{SA}^n$ is convex and $p$ is a positive convex function of $C_{SA}^n$, that is, \;
$\alpha p({\bf x})+(1-\alpha)p({\bf y}) \geq p(\alpha{\bf x}+(1-\alpha){\bf y})>0$, \; for any  ${\bf x}, {\bf y} \in C_{SA}^n$  and $\alpha \in [0,1]$. 
\item[(ii)] 
	If $(x_1,\cdots x_n)\in C_{SA}^n$ and $x_i' \geq x_i$ for all $i=1, \cdots, n$, then  $(x_1',\cdots,x_n')\in C_{SA}^n$.
\end{itemize} 
\end{proposition} 
 
\medskip

\noindent {\it Proof:} (i) We prove by induction on $n$. For $n=1$, $C_{SA}^1=(0,\infty)$ is clearly convex. The convexity and positivity of $p$ amounts to showing that for $x,y>0$, $\alpha \in [0,1]$, 
$${\displaystyle \alpha \frac{1}{x} +(1-\alpha)\frac{1}{y} \geq \frac{1}{\alpha x+(1-\alpha)y}>0 }.$$
This follows from the positivity and convexity of the function $g(x)=1/x$ for $x>0$.
Now suppose that (i) holds for $C_{SA}^{n-1}$  and $W'_{\bf x}={\left(\begin{array}{cc} a' & b' \\ c' & d'\\ \end{array}\right)}=
A_2A_3\cdots A_n$. Let $p'=b'/d'$. 
We again regard the entries  $a'$, $b'$, $c'$, $d'$  and the ratio $p'$ as 
real valued functions of ${\bf x}$. Then since 
$$
A_1 = {\left(\begin{array}{cc} 0 & 1 \\ -1 & x_1\\ \end{array}\right)}, 
$$ 
we have 
$$W_{\bf x}={\left(\begin{array}{cc} c'({\bf x}) & d'({\bf x}) \\ -a'({\bf x})+c'({\bf x}) x_1 & -b'({\bf x})+d'({\bf x}) x_1\\ \end{array}\right)}. \eqno(*)$$
Since $W_{\bf x}$ is strictly admissible, by proposition \ref{Prop:CharacterizeAdm},
$${\displaystyle -b'({\bf x})+d'({\bf x}) x_1>0 \Longrightarrow x_1>\frac{b'({\bf x})}{d'({\bf x})}=p'({\bf x})}.$$ 
Similarly, for  ${\bf y} = (y_1, \cdots, y_n)$  we have 
 $${\displaystyle y_1>\frac{b'({\bf y})}{d'({\bf y})}=p'({\bf y})}.$$ 
Hence 
\begin{eqnarray*}
\alpha x_1+(1-\alpha )y_1 &>& \alpha p'({\bf x})+(1-\alpha )p'({\bf y})\\
&\geq& p'(\alpha {\bf x}+(1-\alpha){\bf y}) \qquad \hbox{(by the hypothesis)}\\
&=& \frac{b'(\alpha {\bf x}+(1-\alpha){\bf y})}{d'(\alpha {\bf x}+(1-\alpha){\bf y})}.\end{eqnarray*}
By proposition \ref{Prop:ByX} (i) this implies that $W_{\alpha {\bf x}+(1-\alpha){\bf y}}$ is strictly admissible, hence $C_{SA}^n$ is convex. Also, proposition \ref{Prop:CharacterizeAdm} implies that $p({\alpha {\bf x}+(1-\alpha){\bf y}})>0$.

Now by $(*)$, 
$${\displaystyle p({\bf x})=\frac{d'({\bf x})}{-b'({\bf x})+d'({\bf x}) x_1}=\frac{1}{-p'({\bf x})+x_1}  },$$
similarly,
$${\displaystyle p({\bf y})=\frac{1}{-p'({\bf y})+y_1}}$$ 
and 
$$p(\alpha {\bf x}+(1-\alpha){\bf y})=\frac{1}{-p'(\alpha {\bf x}+(1-\alpha){\bf y})+\alpha x_1+(1-\alpha)y_1 }.$$
We have \begin{eqnarray*}
&&p(\alpha {\bf x}+(1-\alpha){\bf y})\\
&=&\frac{1}{-p'(\alpha {\bf x}+(1-\alpha){\bf y})+\alpha x_1+(1-\alpha)y_1} \\
&<&\frac{1}{-\alpha p'({\bf x})+(\alpha-1)p'({\bf y})+\alpha x_1+(1-\alpha)y_1}, \\
&& {\vbox to 7mm{}} \hspace{45mm} \hbox{ by induction hypothesis} \\
&=& \frac{1}{{\alpha}\frac{1}{ p({\bf x})}+(1-\alpha)\frac{1}{{p({\bf y}) }}}\\
&\leq& \alpha p({\bf x})+(1-\alpha)p({\bf y}), \qquad \hbox {by convexity of }g(x)=\frac{1}{x}.\\
\end{eqnarray*}
Hence $p$ is a convex function of $C_{SA}^n$ and the proof by induction is complete.

\medskip
(ii) Clearly, it suffices to show that if we increase any entry of a strictly admissible cross ratio vector, the resulting vector remains strictly admissible. This is easily seen by using a suitably normalized configuration of circles for the vector. If the $k$th entry of the vector $(x_1,\cdots,x_n)$ is increased, we use the configuration $\{\overline{C}$, $C_0$, $C_1$, $\cdots$, $C_{n+1}\}$ with $\overline{C}$ the real line, $C_k$ the line $Im(z)=1$ and $C_{k+1}$ the circle $|z-\sqrt{-1}/2|=1/2$. If $p_i$ is the point of tangency of $C_i$ with $\overline{C}$, then we have: 
$p_{k+1}=0$, $p_k=\infty$, $p_{k-1}=x_k$, (by our normalization), 
$p_{k+1}<p_{k+2}<\cdots <p_{n+1}$, \quad $p_{n+1}<p_0<p_1<\cdots<p_{k-1}$, since $(x_1,\cdots,x_n)$ is strictly admissible. The effect of increasing the entry $x_k$ by a positive constant $k$ on the configuration is then to shift  the circles $C_0,C_1,\cdots,C_{k-1}$ to the right by $k$ and to leave all the other circles fixed, if $p_{n+1}<p_0$ to start with, then the inequality remains after this shift and so the configuration remains strictly admissible.
\qed


\subsection{Main lemma}

\medskip
In this subsection, we prove our main lemma by gathering the results stated above.  

Suppose we are given a surface with projective structure  $(\Sigma_g, \mu)$  on 
which a circle packing  $P$  lies.  
Let  $\tau$  be the nerve of the packing  $P$.  
The circle packing  $P$  lifts to a packing 
on the universal cover $\widetilde{\Sigma}_g$, similarly, 
$\tau$ lifts to a triangulation $\widetilde{\tau}$ of 
$\widetilde{\Sigma}_g$.  
For each edge of  $\widetilde{\tau}$, 
the related circles form a configuration on  $\widetilde{\Sigma}_g$. 
Developing them to the Riemann sphere, 
we obtain a configuration of 4 circles on  $\hat{\mathbb C}$  to 
which we assign the cross ratio.  
This assignment is equivariant with respect to 
the action of the covering transformations.  

\begin{definition}[Cross Ratio Parameter and 
Associated Words]\label{Def:CrossRatioParameter}
The map  ${\bf c}_P : E(\tau) \to {\mathbb R}_+$  in the introduction, 
which we call a cross ratio parameter of  $P$,  is the assignment 
defined above. 
Moreover for each vertex  $v \in V(\tau)$  with valence  $n$, 
the associated word  $W_v$  of 
the associated matrices of length  $n$  is obtained 
by reading off the associated matrices of 
the edges incident to   $v$  in the clockwise direction. 
\end{definition} 

Here we have the main lemma. 
\bigskip

\begin{lemma}[Main Lemma]\label{Lem:MainLemma} 
Suppose  ${\bf c}_P$  is a cross ratio parameter of 
some packing  $P$  on a surface with projective structure 
and let  $W_v$  be the associated word for each  $v \in V(\tau)$.  
Then, we have  
\begin{enumerate}
\item[$(1)$] 
	$W_v=-I$, where $I$ is the identity matrix, and 
\item[$(2)$] 
	if the length of  $W_v$  is  $n$, 
	then every subword of $W_v$ of length $\leq n-1$ is admissible and 
	every subword of length  $\leq n-2$  is strictly admissible. 
\end{enumerate}
Conversely, 
if a map  ${\bf c} : E(\tau) \to {\mathbb R}_+$  satisfies  $(1)$  and  $(2)$, 
then there is a unique projective structure on  $\Sigma_g$  together 
with a circle packing  $P$  with nerve isotopic to  $\tau$  
such that the cross ratio parameter  ${\bf c}_P$  is equal to  ${\bf c}$. 
\end{lemma}

\medskip

\noindent {\it Proof:} The proof follows almost directly from propositions \ref{Prop:SurroundConfig} and \ref{Prop:CharacterizeAdm}. We first prove the first half.  
Let $v$ be any vertex of $\tau$ with valence $n$ and let $C$ be the circle corresponding to $v$. The developed image of $C$ is surrounded by $n$ circles $C_1,\cdots, C_n$ (see Figure \ref{Fig:Surround}). Let  $(x_1,\cdots,x_n)$  be the corresponding  cross ratio vector, and $W_v=W_1 \cdots W_n$ the word of associated matrix. By proposition \ref{Prop:SurroundConfig}, $W_v(I_0)=I_0$ where $I_0$ is the interstice between $C_n$, $C$ and $C_1$ when the configuration is in standard position. Hence $W_v=\pm I$.  Furthermore, as the interstices wrap around the central circle $C$ exactly once (see Figure \ref{Fig:Surround}),   any proper subword of $W_v$ is admissible and any subword of length $\leq n-2$ is strictly admissible.  It follows that $W_{n-1}=A_1A_2 \cdots A_{n-1}=-A_n^{-1}$ since $W_{n-1}A_n=\pm I$ and the $(1,2)$ entry of $W_{n-1}$ is negative. Hence $W_v=-I$. 

For the last half, given a map  ${\bf c} : E(\tau) \to {\mathbb R}_+$, the conditions of the lemma are exactly what we need to uniquely define a developing map from $\widetilde{\Sigma}_g$ to $\hat{\mathbb C}$ which is a local homeomorphism and which gives a circle packing on $\widetilde{\Sigma}_g$  up to the composition of elements of  $PSL_2({\mathbb C})$.  This gives the required projective structure and circle packing. \qed

\begin{figure}[ht]
  \begin{center}
    \begin{picture}(140,120)
      \put(0,-20){\scalebox{0.7}{\includegraphics{Surround.eps}}}
\put(0,-20){
\put(66,67){$C$}
\put(66,110){$C_1$}
\put(97,93){$C_2$}
\put(99,65){$C_3$}
\put(91,41){$C_4$}
\put(33,92){$C_n$}
}
    \end{picture}
  \end{center}
  \caption[]{}
\label{Fig:Surround}
\end{figure}

\medskip

We recall  the following definition of the cross ratio parameter 
space mentioned in the introduction. 

\medskip 

\begin{definition}[Cross Ratio 
Parameter Space]\label{Def:CrossRatioParameterSpace}
The set of all maps in ${\mathbb R}^{E(\tau)}$ satisfying 
the conditions (1) and (2) for each vertex of  $\tau$  in the main lemma 
will be called the cross ratio parameter space 
and denoted by
\begin{equation*} 
	{\mathcal C}_{\tau} = 
	\{ {\bf c} : E(\tau) \to {\mathbb R}_+ \, \vert \, 
	{\bf c} \; \; \text{satisfies (1) and (2) for each vertex} \}. 
\end{equation*}  
\end{definition} 
\medskip 

\noindent 
{\bf Remark:} 
It is worth noting that if  $\tau$  is covered by a non-simple 
graph in the universal cover, 
then there are no maps of  $E(\tau)$  to  ${\mathbb R}_+$  satisfying 
the conditions  $(1)$  and  $(2)$  and hence 
${\mathcal C}_{\tau}$  is an empty set. 
\medskip 


\subsection{Holonomy representation}

\medskip

\ \ Given a surface with projective structure together with a  circle packing with nerve $\tau$, we have a cross ratio parameter defined on the edges of $\tau$. An extension of proposition \ref{Prop:SurroundConfig} allows us to find the holonomy representation using the cross ratios and the combinatorics of $\tau$ (compare with the Fenchel-Nielsen coordinates where it is generally quite difficult to obtain the holonomy representation from the coordinates). To begin, we fix the developing map and the holonomy by developing a fixed  interstice to the standard interstice. The holonomy representation of each covering transformation $\gamma \in \pi_1(\Sigma_g)$ is given by a word in the alphabet $\{A_i^{\pm 1}, R^{\pm 1}\}$ where $A_i$ are the associated matrices of the cross ratios and $R={\left(\begin{array}{cc} 0 & \sqrt{-1} \\ \sqrt{-1} & 1\\ \end{array}\right)}\in PSL_2({\mathbb C})$ fixes the standard interstice ${\mathcal I}_s$ and permutes its vertices in an anti-clockwise direction ($R(0)=\sqrt{-1}$, $R(\sqrt{-1})=\infty$, $R(\infty)=0$).
We start with some definitions:

\medskip

\begin{definition} 
A {\it marked triangle} $(\Delta, v)$ of $\widetilde{\Sigma_g}-\tilde{\tau}$ is a triangle with a marked vertex. This corresponds to an interstice of the circle packing on $\widetilde{\Sigma_g}$ with a marked edge. There are two types of basic moves on the set of marked triangles defined as follows:

\smallskip
\begin{itemize} 
\item[(i)]
A type 1 move on the set of marked triangles is a rotation either clockwise or anti-clockwise about the marked vertex of a marked triangle to an adjacent triangle. This moves a marked triangle to an adjacent marked triangle with marked vertex  in the same position. If $x$ is the cross ratio of the edge common to the two triangles, we associate the matrix $A={\left(\begin{array}{cc} 0 & 1 \\ -1 & x\\ 
\end{array}\right)}$ to a  clockwise rotation, and $A^{-1}$ to an anti-clockwise rotation.
\item[(ii)]
A  type 2 move fixes the triangle but changes its marking by rotating the vertices  either anti-clockwise or clockwise. We associate the matrix $R={\left(\begin{array}{cc} 0 & \sqrt{-1} \\ \sqrt{-1} & 1\\ \end{array}\right)}$ for an anti-clockwise rotation of the marking and $R^{-1}$ for a clock-wise rotation.
\end{itemize} 
\end{definition} 

\medskip

Suppose that $\gamma(\Delta_1, v_1)=(\Delta_2, v_2)$ for some $\gamma \in \pi_1(\Sigma_g)$. Then the holonomy representation $\rho(\gamma)$ of $\gamma$  acts on the developing image of  $(\Delta_1, v_1)$ by $\rho(\gamma)(dev(\Delta_1, v_1))=(dev(\Delta_2, v_2))$. Fix a marked triangle $(\Delta_0,v_0)$ and develop it to the standard interstice ${\mathcal I}_s$ such that $v_0$ corresponds to the real line. We have the following result.

\medskip

\begin{proposition}\label{Prop:Holonomy}
Let $m_1, m_2, \cdots, m_{k_i}$ be a sequence of moves that maps the marked triangle $(\Delta_0,v_0)$ to $(\Delta_i, v_i)$ in  $\widetilde{\Sigma_g}-\tilde{\tau}$ and let $B_1,  \cdots$, $B_{k_i}$ be the associated matrices of these moves. Suppose further that $dev$ $(\Delta_0,v_0)=\mathcal {I}_s$ with $v_0$ corresponding to the real line. Then the transformation in $PSL_2(\mathbb{C})$ taking $dev(\Delta_0,v_0)=\mathcal {I}_s$ to $dev(\Delta_i,v_i)$ is $W_i=B_1B_2 \cdots B_{k_i}$. Hence if $\gamma \in \pi_1(\Sigma_g)$ and  $\gamma(\Delta_1, v_1)=(\Delta_2, v_2)$, then the holonomy representation $\rho$ of $\gamma$ is given by
$\rho(\gamma)= W_2W_1^{-1}$.
\end{proposition}

\medskip

The proof of the above is essentially the same as the proof of proposition \ref{Prop:SurroundConfig} and will be omitted.

\medskip

\begin{example}\label{Ex:Holonomy} 
{\em
Consider the triangulation $\tau$ of the torus $\Sigma_1$ with one vertex and three edges, $e_1,e_2$ and $e_3$ and suppose we have a projective structure on the torus which admits a circle packing with nerve $\tau$.   
A fundamental domain for the torus consists of two triangles, see Figure \ref{Fig:Hexagonal} for a portion of the triangulation of $\widetilde{\Sigma_1}$ where all horizontal edges are lifts of $e_1$, all edges of angle $2\pi/3$ with the positive $x$-axis are lifts of $e_2$ and all edges of angle $\pi/3$ with the positive $x$-axis are lifts of $e_3$. Let $x_1,x_2$ and $x_3$ be the cross ratios and $A_1$, $A_2$ and $A_3$ be the associated matrices of $e_1, e_2$ and $e_3$ respectively. The fundamental group of $\Sigma_1$ is generated by three elements $\gamma_1$, $\gamma_2$ and $\gamma_3$  where $\gamma_1(\Delta_0, v_0)=(\Delta_2,v_2)$,
and $\gamma_2(\Delta_1, v_1)=(\Delta_3,v_3)$, $\gamma_3(\Delta_2, v_2)=(\Delta_4,v_4)$, and $\gamma_3\gamma_2^{-1}\gamma_1=id$ (see Figure \ref{Fig:Hexagonal}). It follows from proposition \ref{Prop:Holonomy} that
\begin{eqnarray*}
\rho(\gamma_1)&=& A_3A_1R\\
\rho(\gamma_2)&=& A_3A_1A_2RA_3^{-1}\\
\rho(\gamma_3)&=& A_3A_1A_2A_3RA_1^{-1}A_3^{-1}.\\
\end{eqnarray*}
}
\end{example} 
\begin{figure}[ht]
  \begin{center}
    \begin{picture}(140,135)
      \put(0,-5){\scalebox{0.7}{\includegraphics{Hexagonal.eps}}}
\put(65,86){$\Delta_0$}
\put(84,76){$\Delta_1$}
\put(84,52){$\Delta_2$}
\put(65,43){$\Delta_3$}
\put(45,52){$\Delta_4$}
\put(65,76){$v_0$}
\put(80,69){$v_1$}
\put(84,44){$v_2$}
\put(61,36){$v_3$}
\put(46,43){$v_4$}
    \end{picture}
  \end{center}
  \caption[]{}
\label{Fig:Hexagonal}
\end{figure}

%
%

\section{Around the Andreev-Thurston Solution} 


\subsection{Formal dimension count}

Given a triangulation  $\tau$  on a surface  $\Sigma_g$, 
we easily conclude from the main lemma in 
the previous section that: 

\medskip 

\begin{proposition} 
The cross ratio parameter space ${\mathcal C}_{\tau}$  
is a real semi-algebraic set of formal dimension  $6g-6$  
when  $g \geq 1$.  
\end{proposition} 

\medskip 

\noindent {\it Proof:} 
Since the two conditions in the main lemma are expressed in terms of 
the identities and inequalities in cross ratios, 
 ${\mathcal C}_{\tau}$ is a real semi-algebraic set.
It remains to do a  formal dimension count. 
An easy application of the Gauss-Bonnet theorem shows that 
$|E(\tau)|=6g-6+3|V(\tau)|$. 
For each vertex $v$, the condition $W_v=-I$ from lemma \ref{Lem:MainLemma} 
gives 3 equations in the cross ratios.  
Since there are $|V(\tau)|$ vertices, this gives $3|V(\tau)|$ equations. 
The formal dimension is, by definition, the number of variables 
minus the number of equations and we thus have 
$|E(\tau)|-3|V(\tau)|=6g-6$.
\qed
\medskip 

\noindent
{\bf Remark:} 
The formal dimension count works even for 
$\tau$  which is covered by a non-simple graph in the universal cover, 
because the computation uses only euler characteristic argument.  
However, 
the starting point of the argument in the subsequent subsections 
is the Andreev-Thurston solution, 
so that  $\tau$  will have to be covered by 
a simple graph in the universal cover in the sequel. 
\medskip 


\subsection{When  $ g \geq 2$} 

In this subsection, 
we will see, in the case  $g \geq 2$,  that 
the formal dimension of  ${\mathcal C}_{\tau}$  is the  
correct one at least near the Andreev-Thurston solution and 
that the solution has a nice manifold neighborhood.  
The arguments are based heavily on the well developed 
deformation theory of Kleinian groups by Ahlfors, Bers and 
other authors, see for instance  \cite{Bers}.  

When  $g \geq 2$, the image of the canonical section  $s : {\mathcal T}_g 
\to {\mathcal P}_g$  has a tubular neighborhood 
which consists of projective structures 
whose developing map is injective 
and whose developed image is surrounded by a topological circle
embedded in  $\hat{\mathbb C}$.  
This neighborhood is called the quasi Fuchsian space and 
is denoted by  ${\mathcal Q}{\mathcal F}_g$. 
The projective structure in  ${\mathcal Q}{\mathcal F}_g$  is 
sometimes called quasi Fuchsian. 
The bounding circle of the developed image 
of a quasi Fuchsian structure bounds a disk on the other side in 
$\hat{\mathbb C}$  by the Jordan curve theorem 
and the action of the holonomy image there
is still properly discontinuous. 
Hence by taking the quotient, we obtain another Riemann surface homeomorphic 
to  $\Sigma_g$. 
By assigning to each quasi Fuchsian structure, 
its conformal structure with opposite orientation and 
the conformal structure induced on the other side, 
we obtain a map
\begin{equation*} 
	\beta : {\mathcal Q}{\mathcal F}_g \to 
		\overline{\mathcal T}_g \times {\mathcal T}_g. 
\end{equation*} 
Here we adopt the convention 
so that the original projective structure with opposite orientation
projects down to the first component.  
Bers' simultaneous uniformization theorem states that  
$\beta$  is in fact a homeomorphism. 

\medskip 

\begin{lemma}\label{Lem:NbdAT2}
Suppose  $g \geq 2$, 
and let  $\beta_2$  be the composition of  $\beta$  with 
the projection to the second factor. 
Then  
\begin{equation*} 
	\beta_2 \circ f : f^{-1}({\mathcal Q}{\mathcal F}_g) \to {\mathcal T}_g 
\end{equation*} 
is a homeomorphism, where $f$ is the ``forgetting map'' from ${\mathcal C}_{\tau}$ to ${\mathcal P}_g$ defined in the introduction.  
In particular, 
$f^{-1}({\mathcal Q}{\mathcal F}_g) \subset {\mathcal C}_{\tau}$, 
which includes the Andreev-Thurston 
solution, is homeomorphic to the euclidean space of 
dimension  $6g-6$. 
\end{lemma} 

\medskip 

\noindent
{\it Proof}. 
We first construct a hyperbolic 3-manifold
using the Andreev-Thurston solution.  
Given  $\tau$, the Andreev-Thurston solution gives a 
hyperbolic structure on  $\Sigma_g$  and a packing  $P_0$.  
Developing the hyperbolic surface, 
we obtain a packing on the upper half plane  ${\mathbb H}^2$.  
Let  $R$  be its mirror image on the lower half plane 
$\overline{\mathbb H}^2$.  
The dual packing  $R^*$   which consists of 
the circumscribed circles of interstices in  $R$  fills  
$\overline{\mathbb H}^2$.  
Each member of both packings bounds a hemisphere in 
${\mathbb H}^3$  and such a hemisphere bounds a hemiball  
facing  $\overline{\mathbb H}^2$.  
Chopping off those hemiballs from  ${\mathbb H}^3$, 
we obtain a hyperbolic 3-manifold  $L$  with 
ideal polygonal boundary.  
The point of tangency in the packing  $R$  becomes an ideal vertex 
with rectangler section, 
and its modulus is the cross ratio of a corresponding 
point of tangency in  $P_0$  times  $\sqrt{-1}$. 

The holonomy image of the Andreev-Thurston solution 
acts properly discontinuously on  ${\mathbb H}^3$  by the Poincar\'e 
extension and  $L$  is invariant under its action. 
Hence taking quotient, 
we obtain a hyperbolic 3-manifold  $M$  with ideal 
polygonal boundary and with one end homeomorphic to  $\Sigma_g$.  
The polygonal faces on the boundary are divided
into two classes, 
one extends to the hemisphere bounded by a circle in  $R$  and 
the other in  $R^*$.  
Take the double of  $L$  first along faces in the first class 
and then again take the double of the result of the first doubling  
along faces in the second class. 
This double doubling construction gives us 
a geometrically finite hyperbolic 3-manifold  $N$  with $4$ ends, 
each homeomorphic to  $\Sigma_g$,  and 
$|E(\tau)|$ cusps corresponding to the points of 
tangency in  $P_0$. 
It is obvious by the construction that the modulus of 
each cusp in  $N$  is the cross ratio of the corresponding 
point of tangency in  $P_0$  times  $\sqrt{-1}$. 
Moreover, $N$  naturally admits a  
${\mathbb Z}/2 \times {\mathbb Z}/2 \; (=G)$  symmetry generated 
by two reflections coming from the double doubling.  

We now let  $Def^G(N)$  be the space of 
geometrically finite deformations of  $N$  preserving  the $G$ action.  
Again by the theorem of Bers, 
it can be identified with the space of $G$-invariant quasi conformal 
deformations of the end of  $N$   and 
thus we have a homeomorphism 
\begin{equation*} 
	q : Def^G(N) \to {\mathcal T}_g. 
\end{equation*} 

In any geometrically finite deformation of  $N$,  
the cusps remain being cusps and  
there are  $|E(\tau)|$  cusps.  
Moreover, 
since the deformation is chosen to admit 
a $G$ action generated by two reflections, 
each cusp must admit the same symmetry 
so that the fundamental domain of the cusp remains rectangular.  
Thus for each deformation of  $N$,  
we have a map assigning to each edge of  $\tau$  
a positive real number obtained by the modulus of 
the corresponding cusp over  $\sqrt{-1}$.  
This gives a map  $m : Def^G(N) \to {\mathbb R}^{E(\tau)}$.  

Let  $N'$  be a $G$-invariant geometrically finite 
deformation of  $N$.  
Then the quotient  $M'$  of  $N'$  by the  $G$ action has 
an ideal polygonal boundary having the same combinatorial 
pattern with that of  $M$.  
The universal cover  $L'$  of  $M'$  also 
has an ideal polygonal boundary having the same 
combinatorial pattern with that of  $L$.  
The faces corresponding to  $R$  extend to hemispheres 
which provide a circle packing on the domain 
bounded by a Jordan curve in  $\hat{\mathbb C}$.    
The packing is obviously invariant by the action of 
the holonomy image of  $N'$  and thus by taking quotient 
we obtain a  projective structure on  $\Sigma_g$  together 
with a circle packing with nerve isotopic to  $\tau$. 
The image of a deformation  $N'$  by  $m$  is nothing but 
the cross ratio parameter of the packing so obtained. 
Thus we have squeezed the image of  $m$  by 
\begin{equation*} 
	m : Def^G(N) \to {\mathcal C}_{\tau}. 
\end{equation*} 
Now it is not hard to see that 
the composition of  $m \circ q^{-1}$  and  $b_2 \circ f$, 
is the identity and so is the composition with the reversed order. 
\qed
 

\subsection{When  $g = 1$} 

In this subsection, 
we show that the Andreev-Thurston solution has 
a nice manifold neighborhood of dimension  $2$  when  $g=1$, 
in contrast with the formal dimension count 
of  ${\mathcal C}_{\tau}$  in \S 3.1. 
The arguments are based heavily on the hyperbolic Dehn surgery 
theory developed by Thurston \cite{Thu}. 

The fundamental group of  $\Sigma_1$  is isomorphic 
to  ${\mathbb Z} \times {\mathbb Z}$  and 
is mapped by the holonomy representation of 
a projective structure to an elementary subgroup in 
$PSL_2({\mathbb C})$  which fixes a single point or two points 
on  $\hat{\mathbb C}$  globally.  
Thus since  ${\mathcal C}_{\tau}$  is the space of 
deformations of projective structures on  $\Sigma_1$  with 
a specific type of packings, 
we in particular obtain a family of such representations of 
${\mathbb Z} \times {\mathbb Z}$  in  $PSL_2({\mathbb C})$  
parameterized by  ${\mathcal C}_{\tau}$.  
In this setting, 
the hyperbolic Dehn surgery theory defines a continuous map 
\begin{equation*} 
	d : {\mathcal C}_{\tau} \to {\mathbb R}^2 \cup \{ \infty \}/\sim 
\end{equation*} 
which generalizes the classical Dehn surgery 
coefficients in the knot theory, 
such that 
the Andreev-Thurston solution is mapped to  $\infty$  by  $d$. 
The map  $d$  is in fact described in terms of the complex affine 
structure in \cite{Thu} and the image lies in 
${\mathbb R}^2 \cup \{ \infty \}$, 
but  ${\mathcal C}_{\tau}$  here is described in terms of 
the projective structure which is generically 
in one to two correspondence with the complex affine structure, 
so that we divide the image by 
the antipodal involution  $\sim$  of  ${\mathbb R}^2$. 

\medskip 

\begin{lemma}\label{Lem:NbdAT1}
Suppose that  $g =1$.  
Then there is a neighborhood  $U$  of the Andreev-Thurston solution 
in  ${\mathcal C}_{\tau}$  which is homeomorphic to 
the euclidean space of dimension  $2$,  
and the restriction of  $d$  to which is an embedding.  
\end{lemma} 

\medskip 

\noindent
{\it Proof}. 
Given  $\tau$, the Andreev-Thurston solution gives an   
euclidean torus and a packing  $P_0$.  
Developing the torus, 
we obtain a packing  $R$  on the complex plane  ${\mathbb C}$.  
Using the same dualizing, chopping off, taking quotients and double doubling 
construction in the previous subsection, 
we obtain a dual packing  $R^*$  on  ${\mathbb C}$  and 
the hyperbolic 3-manifolds  $L, \; M, \; N$  respectively, 
where  $N$  has  $4$  old cusps  and 
$|E(\tau)|$ new cusps corresponding to the points of 
tangency in  $P_0$.  
It is again obvious by the construction that the modulus of 
each cusp is the cross ratio of the corresponding 
point of tangency in  $P_0$  times  $\sqrt{-1}$. 
Moreover, $N$  admits a  
${\mathbb Z}/2 \times {\mathbb Z}/2 \; (=G)$  symmetry generated 
by the two reflections coming from the double doubling.  

Let  $Def^G(N)$  be in this case a small neighborhood of 
$N$  in the space of hyperbolic Dehn surgery deformations 
of  $N$  preserving  the $G$ action.  
Then since the deformations preserve the $G$ action, 
the cusps created by the double doubling 
must admit the  $G$-symmetry generated by the two reflections, 
so they remain cusps.  
Thus by Thurston's theorem, 
the deformations of  $N$  which are small enough are parameterized by the 
$G$-invariant generalized Dehn surgery coefficients 
of the other old cusps and we have an embedding 
\begin{equation*} 
	q : Def^G(N) \to {\mathbb R}^2 \cup \{ \infty \}/\sim 
\end{equation*} 
where the image is a neighborhood of  $\infty = q(N)$. 

In any Dehn surgery deformation of  $N$  in  $Def^G(N)$,  
there are  $|E(\tau)|$  cusps corresponding to the points 
of tangency in  $P_0$. 
Moreover since such cusps admit reflectional symmetry, 
their fundamental domains remain rectangular.  
Thus for each deformation of  $N$,  
we have a map assigning to each edge of  $\tau$  a positive 
real number obtained by the modulus of 
the corresponding cusp over  $\sqrt{-1}$.  
Then the assignment of this map to  $N$  defines 
a map  $m : Def^G(N) \to {\mathbb R}^{E(\tau)}$.  
 
Let  $N'$  be a $G$-invariant incomplete Dehn surgery deformation 
of  $N$  whose completion lies in  $Def^G(N)$.  
Thus  $N'$  is still homeomorphic to  $N$.  
Then the quotient  $M'$  of  $N'$  by the  $G$ action has 
an ideal polygonal boundary having the same combinatorial 
pattern with that of  $M$.  
The universal cover  $L'$  of  $M'$  also 
has an ideal polygonal boundary having the same 
combinatorial pattern with that of  $L$.  
Since  $L'$  carries a hyperbolic metric though incomplete, 
it can be developed to  ${\mathbb H}^3$.  
The developed image of a face corresponding to  $R$  extends to 
a hemisphere whose boundary defines a circle on   $\hat{\mathbb C}$.  
The coherent collection of such circles provide 
a developed image of a certain circle packing on  $\Sigma_g$.  
Pulling back the projective structure on  $\hat{\mathbb C}$   
to the universal cover of  $\Sigma_g$, 
and taking quotient by the covering transformations, 
we obtain a  projective structure on  $\Sigma_1$  together 
with a circle packing with nerve isotopic to  $\tau$. 
The image of a deformation  $N'$  by  $m$  is nothing but 
the cross ratio parameter of the packing so obtained. 
We  have thus again squeezed the image of  $m$  by 
\begin{equation*} 
	m : Def^G(N) \to {\mathcal C}_{\tau}. 
\end{equation*} 
Now it is not hard to see that 
$m \circ q ^{-1} $  is the inverse of  $d$  in 
a small neighborhood of  $\infty$.  
\qed
\medskip 

This lemma together with Lemma \ref{Lem:NbdAT2}  gives
the second lemma in the introduction. 

%
%

\section{Topology of  ${\mathcal C}_{\tau}$  for One Circle Packings}


\subsection{Combinatorics}

\medskip
For this section and the next, we assume that $g \geq 2$, and $\tau$ is a graph on $\Sigma_g$ with one vertex which triangulates $\Sigma_g$. By the Gauss-Bonnet, $\tau$ has $6g-3$ edges and triangulates $\Sigma_g$ into $4g-2$ triangles. Let $\tau'$ be the dual graph of $\tau$.  $\tau'$ is a trivalent graph (all vertices have valence 3) with $6g-3$ edges $e_1', \cdots, e_{6g-3}'$. The cross ratios defined on the edges of $\tau$ clearly define a map of  $E(\tau')$  to  ${\mathbb R}_+$  which we also call a cross ratio parameter, hence the cross ratio parameter space  ${\mathcal C}_{\tau}$  can also be interpreted as being defined on the edges of $\tau'$, and we will adopt this point of view whenever appropriate.  

Cutting $\Sigma_g$ along the edges of $\tau'$ gives a polygon ${\mathcal P}$ with $12g-6$ sides $s_1, s_2, \cdots s_{12g-6}$ and $6g-3$ side pairings $s_i \leftrightarrow s_j$, each paired side corresponding to some edge $e_s'$ of $\tau'$. When $g=2$ there are essentially $8$ different possible side-pairing patterns, see Figure \ref{Fig:PairingPattern}. When $g=3$, a computer search produces over nine hundred different side-pairing patterns. 

\begin{figure}[ht]
  \begin{center}
    \begin{picture}(280,110)
\put(0,-10){
      \put(0,60){\scalebox{0.3}{\includegraphics{PairingPattern-1.eps}}}
      \put(80,60){\scalebox{0.3}{\includegraphics{PairingPattern-2.eps}}}
      \put(160,60){\scalebox{0.3}{\includegraphics{PairingPattern-3.eps}}}
      \put(240,60){\scalebox{0.3}{\includegraphics{PairingPattern-4.eps}}}
      \put(0,0){\scalebox{0.3}{\includegraphics{PairingPattern-5.eps}}}
      \put(80,0){\scalebox{0.3}{\includegraphics{PairingPattern-6.eps}}}
      \put(160,0){\scalebox{0.3}{\includegraphics{PairingPattern-7.eps}}}
      \put(240,0){\scalebox{0.3}{\includegraphics{PairingPattern-8.eps}}}
}
    \end{picture}
  \end{center}
  \caption[]{}
\label{Fig:PairingPattern}
\end{figure}

Since each vertex of $\tau'$ has valence 3, it lifts to 3 vertices of ${\mathcal P}$ and the three edges incident to the vertex corresponds to $3$ pairs of sides of ${\mathcal P}$. It follows that if $s_i \leftrightarrow s_j$ and $s_{i+1} \leftrightarrow s_k$, then $s_{j-1} \leftrightarrow s_{k+1}$ (indices are taken mod $12g-6$). $\{s_i \leftrightarrow s_j, s_{i+1} \leftrightarrow s_k, s_{j-1} \leftrightarrow s_{k+1}\}$ is called the triple of side-pairings associated to the vertex. We can divide the triples to 2 types according to whether the pairs of sides separate each other as follows:

\medskip

\begin{definition}\label{Def:Triple} 
The triple of side pairings $\{s_i \leftrightarrow s_j, s_{i+1} \leftrightarrow s_k, s_{j-1} \leftrightarrow s_{k+1}\}$ is {\it separating} if $i,i+1,j-1,j,k,k+1$ is in cyclic order mod $12g-6$ and {\it non-separating} if $i,i+1,k,k+1,j-1,j$ is in cyclic order mod $12g-6$ (see Figures \ref{Fig:Triple}a and \ref{Fig:Triple}b). Similarly, a vertex of $\tau'$ and the triple of  edges incident to the vertex are called  separating and  non-separating if they give rise to separating and non-separating triples of side-pairings respectively. 
\end{definition} 

\begin{figure}[ht]
  \begin{center}
    \begin{picture}(200,100)
      \put(0,10){\scalebox{2}{\includegraphics{Triple-a.eps}}}
      \put(120,10){\scalebox{2}{\includegraphics{Triple-b.eps}}}
      \put(0,-4){(a) separating}
      \put(28,89){$i$}
      \put(40,89){$i+1$}
      \put(76,35){$j-1$}
      \put(70,24){$j$}
      \put(3,24){$k$}
      \put(-22,38){$k+1$}

      \put(110,-4){(b) non-separating}
      \put(148,89){$i$}
      \put(160,89){$i+1$}
      \put(196,35){$k$}
      \put(190,24){$k+1$}
      \put(105,24){$j-1$}
      \put(115,38){$j$}
    \end{picture}
  \end{center}
  \caption[]{}
\label{Fig:Triple}
\end{figure}

\medskip
\noindent {\bf Remark:} In the case of a non-separating triple, a simple topological argument shows that $s_{i+1}$ and $s_k$ are not adjacent sides of ${\mathcal{P}}$,  similarly  $s_{k+1}$ and $s_{j-1}$ are not adjacent, and $s_j$ and $s_i$ are not adjacent. 

\medskip

 In the case of genus 1, there is only one triple and it is separating. It will be crucial in our analysis of the one circle packing case for higher genus that there exists a non-separating triple. We have the following result: 

\medskip

\begin{proposition}\label{Prop:NonseparatingExist} 
If $g \geq 2$ and $\tau'$ is a trivalent graph on $\Sigma_g$ cutting $\Sigma_g$ into one polygon, then both separating and non-separating  vertices of $\tau'$ must occur.
\end{proposition}

\medskip

\noindent {\it Proof:} 
We first show that there exists non-separating triples of side-pairings.
If $s_1 \leftrightarrow s_i$ and $s_2 \leftrightarrow s_j$ with $i>j$, then 
$\{s_1 \leftrightarrow s_i, s_2 \leftrightarrow s_j, s_{j+1} \leftrightarrow s_{i-1}\}$ is a non-separating triple and we are done. So we may assume that $i<j$ and  $\{s_1 \leftrightarrow s_i, s_2 \leftrightarrow s_j, s_{j+1} \leftrightarrow s_{i-1}\}$ is a separating triple. Furthermore, by a cyclic permutation if necessary, we may assume that $i \geq 5$ since $6g-3 >6$. Now consider the two side-pairings $s_2 \leftrightarrow s_j$ and
$s_3 \leftrightarrow s_k$. Again, if $k<j$ this gives rise to a non-separating triple $\{s_2 \leftrightarrow s_j, s_3 \leftrightarrow s_k, s_{k+1} \leftrightarrow s_{j-1}\}$, so we may assume again that $k>j$. Continuing inductively, we either get a non-separating triple or we eventually have $s_{i-2}\leftrightarrow s_l$ where $l>j+2$. This gives us the non-separating triple $\{s_{i-2} \leftrightarrow s_l, s_{i-1} \leftrightarrow s_{j+1}, s_{j+2} \leftrightarrow s_{l-1}\}$ and we are done.

To show that there exists a separating triple, We use a similar method. We may assume that
$\{s_1 \leftrightarrow s_i, s_2 \leftrightarrow s_j, s_{j+1} \leftrightarrow s_{i-1}\}$ is a non-separating triple so that $j<i$ and furthermore $j>3$ by the remark following definition \ref{Def:Triple}. Now if $s_3 \leftrightarrow s_k$, either $\{s_2 \leftrightarrow s_j, s_3 \leftrightarrow s_k, s_{k+1} \leftrightarrow s_{j-1}\}$ is a separating triple or $k<j$. After less than $j-1$ such steps, we either get a separating triple or we reach a stage where $s_m \leftrightarrow s_p$, $s_{m+1} \leftrightarrow s_q$ with $q>p$. This gives the separating triple $\{s_{m} \leftrightarrow s_p, s_{m+1} \leftrightarrow s_{q}, s_{p-1} \leftrightarrow s_{q+1}\}$ as required. \qed


\subsection{Free cross ratios}

\medskip

In this subsection, we will show that we can choose a set of $6g-6$ independent cross ratios on the set of edges of $\tau$, (equivalently of $\tau'$) such that they determine the remaining three cross ratios. The idea is to choose a non-separating vertex $ v'$ of $\tau'$ (always possible by proposition \ref{Prop:NonseparatingExist}) and use the cross ratios of the 3 edges of $\tau'$ incident to $v'$  as the set of dependent variables and the remaining $6g-6$ variables as the free variables.
 
We start with a local argument. Recall from \S 2 that apart from the admissibility conditions which are open conditions, locally, the cross ratio parameter space ${\mathcal C}_{\tau}$ is given by the matrix identity $W_v:={\left(\begin{array}{cc} a & b \\ c & d\\ \end{array}\right)}
=-I$ or equivalently, by the 3 equations
\begin{eqnarray*}
a(x_1,\cdots,x_{6g-3})&=&-1\\
b(x_1,\cdots,x_{6g-3})&=&0\\
c(x_1,\cdots,x_{6g-3})&=&0.\\
\end{eqnarray*}
(Note that the entries $a,b,c$ and $d$ are regarded as functions of the variables and that $d=-1$ follows from the above since $\det(W_v)=1$). 
We will show that we can find a set of 3 variables $x_i,x_j$ and $x_k$ such that the Jacobian 
$${\displaystyle {\left|\begin{array}{ccc} \frac{\partial a}{\partial x_i} & \frac{\partial a}{\partial x_j} & \frac{\partial a}{\partial x_k}   \\ 
\frac{\partial b}{\partial x_i} &  \frac{\partial b}{\partial x_j} & \frac{\partial b}{\partial x_k} \\
\frac{\partial c}{\partial x_i} & \frac{\partial c}{\partial x_j} & \frac{\partial c}{\partial x_k}   
\\
\end{array}\right|}} \neq 0.$$  By the implicit function theorem, this implies that locally, ${\mathcal C}_{\tau}$  is determined by the remaining $6g-6$ variables and hence is locally of dimension $6g-6$.

Recall that if $v$ is the vertex of $\tau$ and $e_1, \cdots, e_{6g-3}$ the set of edges with cross ratios $x_1,\cdots, x_{6g-3}$, reading off the edges incident to $v$ in a clockwise direction, we get the cross ratio vector ${\bf x}=(x_{n_1}, x_{n_2}, \cdots,$ $x_{n_{12g-6}})$ and the associated word $W_{\bf x}=A_{n_1}A_{n_2} \cdots A_{n_{12g-6}}$. Alternatively, we can think of the cross ratio vector as the the vector obtained by reading off the cross ratios of the sides of the polygon obtained by cutting $\Sigma_g$ along the dual graph $\tau'$.
\medskip

\begin{example} 
For the polygon with side pairing given in the first pattern of Figure \ref{Fig:PairingPattern},
the cross ratio vector is given by $${\bf x}=(x_1,x_2,x_3,x_4,x_2,x_5,x_3,x_4,x_5,x_1,x_6,x_7,x_8,x_6,x_9,x_7,x_8,x_9).$$ 
\end{example} 

\medskip

By proposition \ref{Prop:NonseparatingExist}, and by doing a cyclic permutation if necessary, we may assume that $\{s_{12g-6} \leftrightarrow s_i, s_1 \leftrightarrow s_j, s_{j+1} \leftrightarrow s_{i-1}\}$ is a non-separating triple so that $j<i$. Let $x, y$ and $z$ be the cross ratios of the edges of $\tau'$ corresponding to the pairs of sides $s_1 \leftrightarrow s_j$, $s_{j+1} \leftrightarrow s_{i-1}$ and $s_{12g-6} \leftrightarrow s_i$ respectively and let $x_1,x_2, \cdots, x_{6g-6}$ be the cross ratios of the remaining edges of $\tau'$ corresponding to the other remaining pairs of sides numbered in some fixed order. The cross ratio vector ${\bf x}$ is then of the form 
$${\bf x}=(x,x_{n_2},\cdots, x_{n_{j-1}},x,y,x_{n_{j+2}},\cdots,x_{n_{i-2}},y,z,x_{n_{i+1}},\cdots,x_{n_{12g-7}},z)$$ 
with associated word $W$ in the associated matrices given by 
\begin{eqnarray*}
&&\!\!\!\!\!\!\!\!\!\!
W_{\bf x}={\left(\begin{array}{cc} a & b \\ c & d\\ \end{array}\right)} \\
&&=
{\tiny {\left(\begin{array}{cc} 0 & 1 \\ -1 & x\\ \end{array}\right)} 
T
{\left(\begin{array}{cc} 0 & 1 \\ -1 & x\\ \end{array}\right)}
{\left(\begin{array}{cc} 0 & 1 \\ -1 & y\\ \end{array}\right)}
U
{\left(\begin{array}{cc} 0 & 1 \\ -1 & y\\ \end{array}\right)}
{\left(\begin{array}{cc} 0 & 1 \\ -1 & z\\ \end{array}\right)}
V
{\left(\begin{array}{cc} 0 & 1 \\ -1 & z\\ \end{array}\right)}}
\end{eqnarray*}
where $T$, $U$ and $V$ are subwords of $W$ of lengths $j-2$, $i-j-3$ and $12g-i-7$ respectively which depend only on the variables  $x_1,\cdots,x_{6g-6}$.

\smallskip

\begin{proposition}\label{Prop:Derivative}
With the notation above, 
let $W_k:={\left(\begin{array}{cc} a_k & b_k \\ c_k & d_k\\ \end{array}\right)}=A_{n_1}A_{n_2}\cdots A_{n_k}$ be the subword of $W_{\bf x}$ given by the product of the first $k$ associated matrices.
Then  \begin{eqnarray*}
{\left(\begin{array}{cc} a_x & b_x \\ c_x & d_x\\ \end{array}\right)}&=&{\left(\begin{array}{cc} 0 & 0 \\ -1 & 0\\ \end{array}\right)}+
{\left(\begin{array}{cc} -b_{j-1}d_{j-1} & b_{j-1}^2 \\ -d_{j-1}^2 & b_{j-1}d_{j-1}\\ \end{array}\right)}\\
{\left(\begin{array}{cc} a_y & b_y \\ c_y & d_y\\ \end{array}\right)}&=&{\left(\begin{array}{cc} -b_{j}d_{j} & b_{j}^2 \\ -d_{j}^2 & b_{j}d_{j}\\ \end{array}\right)}+{\left(\begin{array}{cc} -b_{i-2}d_{i-2} & b_{i-2}^2 \\ -d_{i-2}^2 & b_{i-2}d_{i-2}\\ \end{array}\right)}\\
{\left(\begin{array}{cc} a_z & b_z \\ c_z & d_z\\ \end{array}\right)}
&=&{\left(\begin{array}{cc} -b_{i-1}d_{i-1} & b_{i-1}^2 \\ -d_{i-1}^2 & b_{i-1}d_{i-1}\\ \end{array}\right)}+{\left(\begin{array}{cc} 0 & 1 \\ 0 & 0\\ \end{array}\right)}\\
\end{eqnarray*}
where  $a_x, a_y, a_z$ etc. denotes the partial derivatives of the entries of $W$ with respect to the variables $x$, $y$ and $z$.
\end{proposition} 

\medskip

\noindent {\it Proof:} We will show the second formula, the first and third are similar. We apply  the product rule,  and using the fact that $W_{\bf x}=-I$ since the cross ratio vector corresponds to a point in ${\mathcal C}_{\tau}$, we get
\begin{eqnarray*}
{\left(\begin{array}{cc} a_y & b_y \\ c_y & d_y\\ \end{array}\right)}&=&
W_{j}
{\left(\begin{array}{cc} 0 & 0 \\ 0 & 1\\ \end{array}\right)}
(-W_{j+1})^{-1}+
W_{i-2}
{\left(\begin{array}{cc} 0 & 0 \\ 0 & 1\\ \end{array}\right)}
(-W_{i-1})^{-1}
\\
&=& {\left(\begin{array}{cc} a_j & b_j \\ c_j & d_j\\ \end{array}\right)}
{\left(\begin{array}{cc} 0 & 0 \\ 0 & 1\\ \end{array}\right)}
{\left(\begin{array}{cc} -d_{j+1} & b_{j+1} \\ c_{j+1} & -a_{j+1}\\ \end{array}\right)}\\
&{}& \quad +
{\left(\begin{array}{cc} a_{i-2} & b_{i-2} \\ c_{i-2} & d_{i-2}\\ \end{array}\right)}
{\left(\begin{array}{cc} 0 & 0 \\ 0 & 1\\ \end{array}\right)}
{\left(\begin{array}{cc} -d_{i-1} & b_{i-1} \\ c_{i-1} & -a_{i-1}\\ \end{array}\right)}\\
&=&{\left(\begin{array}{cc} b_jc_{j+1} & -b_ja_{j+1} \\ d_jc_{j+1} & -d_ja_{j+1}\\ \end{array}\right)}
+{\left(\begin{array}{cc} b_{i-2}c_{i-1} & -b_{i-2}a_{i-1} \\ d_{i-2}c_{i-1} & -d_{i-2}a_{i-1}\\ \end{array}\right)}\\
&=&{\left(\begin{array}{cc} -b_{j}d_{j} & b_{j}^2 \\ -d_{j}^2 & b_{j}d_{j}\\ \end{array}\right)}+{\left(\begin{array}{cc} -b_{i-2}d_{i-2} & b_{i-2}^2 \\ -d_{i-2}^2 & b_{i-2}d_{i-2}\\ \end{array}\right)}\\
\end{eqnarray*}
since 
$${\left(\begin{array}{cc} a_{j+1} & b_{j+1} \\ c_{j+1} & d_{j+1}\\ \end{array}\right)}={\left(\begin{array}{cc} a_{j} & b_{j} \\ c_{j} & d_{j}\\ \end{array}\right)}{\left(\begin{array}{cc} 0 & 1 \\ -1 & x_{n_{j+1}}\\ \end{array}\right)}$$
so that $c_{j+1}=-d_j$, $a_{j+1}=-b_j$ and similarly $c_{i-1}=-d_{i-2}$, $a_{i-1}=-b_{i-2}$. \qed

\medskip

We are  now ready to prove the local result for the cross ratio parameter space  ${\mathcal C}_{\tau}$.

\medskip

\begin{lemma}\label{Lem:Manifold}
The cross ratio parameter space ${\mathcal C}_{\tau}$ for a triangulation $\tau$ of $\Sigma_g$ ($g \geq 2$) with one vertex is locally homeomorphic to ${\mathbb R}^{6g-6}$. Furthermore, the local parameters can be taken to be the cross ratios on a set of $6g-6$ edges of $\tau'$ such that the remaining $3$ edges are incident to a non-separating vertex of $\tau'$.
\end{lemma} 

\medskip

\noindent {\it Proof:} Adopting the notation of proposition \ref{Prop:Derivative}, it suffices to prove that the Jacobian
$${\displaystyle {\left|\begin{array}{ccc} \frac{\partial a}{\partial x} & \frac{\partial a}{\partial y} & \frac{\partial a}{\partial z}   \\ 
\frac{\partial b}{\partial x} &  \frac{\partial b}{\partial y} & \frac{\partial b}{\partial z} \\
\frac{\partial c}{\partial x} & \frac{\partial c}{\partial y} & \frac{\partial c}{\partial z}   
\\
\end{array}\right|}} \neq 0$$ 
where $x$, $y$ and $z$ are the parameters of the edges incident to a non-separating vertex. By proposition \ref{Prop:Derivative}, this is equivalent to showing that 
$${\displaystyle {\left|\begin{array}{ccc}-b_{j-1}d_{j-1}& -b_jd_j-b_{i-2}d_{i-2}& -b_{i-1}d_{i-1}     \\ 
 b_{j-1}^2& b_j^2+b_{i-2}^2   & b_{i-1}^2+1 \\
 -1-d_{j-1}^2& -d_j^2-d_{i-2}^2 & -d_{i-1}^2   
\\
\end{array}\right|}} \neq 0. $$ 
We will show that the $3$ columns $\phi_1$, $\phi_2$ and $\phi_3$ of the above Jacobian are linearly independent.
By proposition \ref{Prop:CharacterizeAdm}, 
$$b_{j-1},d_{j-1},b_j,d_j,b_{i-2},d_{i-2},b_{i-1},d_{i-1}>0, ~\hbox{and}$$ 
$$p_{j-1}=b_{j-1}/d_{j-1}<p_j=b_{j}/d_{j}<p_{i-2}=b_{i-2}/d_{i-2}<p_{i-1}=b_{i-1}/d_{1-1}$$ since the subwords  $W_k$ of $W_{\bf x}$  are all strictly admissible for $k \leq 12g-4$. 
Let
$$\vec{v}_1={\left(\begin{array}{c} 0 \\ 0\\ -1\end{array}\right)},
\qquad \vec{v}_2={\left(\begin{array}{c} -p_{j-1}\\ p_{j-1}^2\\ -1\end{array}\right)},
\qquad \vec{v}_3={\left(\begin{array}{c} -p_{j}\\ p_{j}^2\\ -1\end{array}\right)},$$
$$\vec{v}_4={\left(\begin{array}{c} -p_{i-2}\\ p_{i-2}^2\\ -1\end{array}\right)},\qquad
\vec{v}_5={\left(\begin{array}{c} -p_{i-1}\\ p_{i-1}^2\\ -1\end{array}\right)},
\qquad \vec{v}_6={\left(\begin{array}{c} 0\\ 1\\ 0\end{array}\right)}. 
$$
We have 
$\phi_1=\vec{v}_1+d_{j-1}^2 \vec{v}_2, \quad
\phi_2=d_{j}^2
\vec{v}_3+
d_{i-2}^2
\vec{v}_4
, \quad \phi_3=d_{i-1}^2
\vec{v}_5+
\vec{v}_6
.$

Note that the vectors $\vec{v}_1, \vec{v}_2, \vec{v}_3, \vec{v}_4, \vec{v}_5$, and  $\vec{v}_6$ all lie on the circular cone passing through the origin and the parabola $\{y=x^2, \; z=-1\}$ and that furthermore they are in anti-clockwise order. Since $\phi_1$ is a positive linear combination of  $\vec{v}_1$ and $\vec{v}_2$, its intersection with the horizontal plane $\{z=-1\}$ lies on the chord joining $(0,0,-1)$ and $(-p_{j-1},p_{j-1}^2,-1)$. Similarly, the intersection of $\phi_2$ with the horizontal plane $\{z=-1\}$ lies on the chord joining $(-p_j,p_j^2,-1)$ and $(-p_{i-2},p_{i-2}^2,-1)$ and the intersection of $\phi_3$ with the horizontal plane $\{z=-1\}$ lies on the half infinite line segment $\{(x,y,z) ~|~x=-p_{i-1}$, $\; y \geq p_{i-1}^2, \; z=-1\}$ . The convexity of the curve $y=x^2$ and the arrangement of the points now implies that the three vectors are linearly independent. 
\qed

\medskip

\noindent {\bf Remark:} The proofs of proposition \ref{Prop:Derivative} and lemma \ref{Lem:Manifold} show that  more general results hold. For example, for lemma \ref{Lem:Manifold}, we need not choose a triple of edges corresponding to a non-separating vertex to be the ones with the dependent variables. In fact, we can use any 3 edges such that the three sets of paired sides in the fundamental polygon $\mathcal P$ corresponding to these 3 edges  are mutually  non-separating. Indeed, one sees how the proof breaks down if the non-separating property is not present, which is why the result does not hold for the genus 1 case. Similarly, for more general triangulations of $\Sigma_g$, if a vertex with valence $n$ is chosen such that no edge is incident to it more than once, then the 3 equations corresponding to that vertex are always independent and we can choose the cross ratios on any set of $n-3$ edges as the free variables for these equations.


\subsection{Topology of  ${\mathcal C}_{\tau}$}

Lemma \ref{Lem:Manifold} shows that ${\mathcal C}_{\tau}$ is a manifold of dimension $6g-6$ but it does not give information about the topology of ${\mathcal C}_{\tau}$. Before proving a global version of lemma \ref{Lem:Manifold}, we first state and prove the following preliminary result:

\medskip

\begin{proposition}\label{Prop:MatrixIdentity} 
Let $A$, $B$ and $C$ be three strictly admissible words of associated matrices. Then
\begin{eqnarray*}
&&ABC=-I
\\
&\Leftrightarrow& \textrm{(2,2) entries of }
\left\{
 \begin{array}{l}
  AB=-C^{-1}
  \\
  BC=-A^{-1}
  \\
  CA=-B^{-1}
 \end{array}
\right.
.
\end{eqnarray*}
Furthermore, if the above holds, then ABC also satisfies condition $(2)$ of the main lemma, that is, every proper subword of ABC is admissible.
\end{proposition}

\medskip

\noindent {\it Proof:}
$(\Rightarrow)$ is clear.

$(\Leftarrow)$ Let $A={\left(\begin{array}{cc} a_1 & a_2 \\ a_3 & a_4\\ \end{array}\right)}$, $B={\left(\begin{array}{cc} b_1 & b_2 \\ b_3 & b_4\\ \end{array}\right)}$ and $C={\left(\begin{array}{cc} c_1 & c_2 \\ c_3 & c_4\\ \end{array}\right)}$.  By proposition \ref{Prop:CharacterizeAdm}, since $A$, $B$ and $C$ are strictly admissible words, $a_1,b_1,c_1 \leq 0$, $a_2,b_2,c_2 >0$, $a_3,b_3,c_3 <0$ and $a_4,b_4,c_4>0$.
Then the conditions on the $(2,2)$ entries implies
\begin{eqnarray}
-c_1&=& a_3b_2+a_4b_4 \\ 
-a_1&=&b_3c_2+b_4c_4 \\ 
-b_1&=& c_3a_2+c_4a_4. 
\end{eqnarray}
It suffices to show that
\begin{eqnarray*}
-c_4&=&a_1b_1+a_2b_3\\
c_2&=&a_1b_2+a_2b_4\\
c_3&=&a_3b_1+a_4b_3
\end{eqnarray*}
as together with $(1)$ this implies $AB=-C$, hence $ABC=-I$.

From $\det C=1$, we have
\begin{eqnarray}
c_2c_3-c_1c_4+1 &=& 0 \nonumber\\
 \qquad \Longrightarrow \quad a_2b_3c_2c_3-a_2b_3c_1c_4+a_2b_3 &=& 0. \hspace{4cm} 
\end{eqnarray}
From $(2)$ and $(3)$,
\begin{eqnarray*}
-b_3c_2&=&a_1+b_4c_4\\
-a_2c_3&=&b_1+a_4c_4\\
\Longrightarrow \quad a_2b_3c_2c_3&=&a_1b_1+(a_1a_4+b_1b_4)c_4+a_4b_4c_4^2.
\end{eqnarray*}
Substituting into $(4)$, we get
$$a_4b_4c_4^2+(a_1a_4+b_1b_4-a_2b_3c_1)c_4+(a_1b_1+a_2b_3)=0. \eqno(5)$$
We claim that $$a_1a_4+b_1b_4-a_2b_3c_1=1+a_4b_4(a_1b_1+a_2b_3). \eqno(6)$$
Substituting  $-c_1=a_3b_2+a_4b_4$ from $(1)$, this is equivalent to
\begin{eqnarray*}a_1a_4+b_1b_4+a_2b_3a_3b_2&=&1+a_4b_4a_1b_1\\ 
 \Longleftrightarrow \qquad a_2a_3b_2b_3&=&(1-a_1a_4)(1-b_1b_4).
\end{eqnarray*}
This is obviously true since
$$a_1a_4-1=a_2a_3, \qquad b_1b_4-1=b_2b_3.$$
Now substituting $(6)$ into $(5)$, we get 
\begin{eqnarray*}
(c_4+(a_1b_1+a_2b_3))(a_4b_4c_4+1)&=&0\\
\Longrightarrow \qquad -c_4&=&a_1b_1+a_2b_3
\end{eqnarray*}
since $a_4b_4c_4+1>0$. \smallskip

Now from $(2)$, we get
\begin{eqnarray*}
b_3c_2 &=& -a_1-b_4c_4\\
&=& -a_1 +b_4(a_1b_1+a_2b_3)\\
&=& a_1(-1+b_1b_4)+a_2b_4b_3\\
&=& a_1b_2b_3+a_2b_4b_3\\
\Longrightarrow \quad c_2 &=& a_1b_2+a_2b_4, \quad \hbox{since}~~b_3<0.
\end{eqnarray*}
Similarly, from $(3)$, we get
\begin{eqnarray*}
a_2c_3 &=& -b_1-a_4c_4\\
&=& -b_1 +a_4(a_1b_1+a_2b_3)\\
&=& b_1(-1+a_1a_4)+a_2a_4b_3\\
&=& b_1a_2a_3+a_2a_4b_3\\
\Longrightarrow \quad c_3 &=& a_3b_1+a_4b_3, \qquad \hbox{since}~~a_2 \neq 0.
\end{eqnarray*}
This completes the first part of the proposition. 
Finally, we note that since $A$, $B$ and $C$ are strictly admissible, if we consider the configuration of circles associated to $A$, $B$ and $C$, the interstices do not go around the central circle completely. Hence for $ABC$, the interstices go around the central circle less than 3 times. However, since $ABC=-I$, the interstices corresponding to $ABC$ go around the central circle exactly an odd number of times. As it is strictly less than $3$,   they go around exactly once and $ABC$ is a word that satisfies both conditions $(1)$ and $(2)$ of the main lemma. 
\qed

\medskip

We now prove the global version of lemma \ref{Lem:Manifold}:

\medskip

\begin{lemma}\label{Lem:Cell}
Let $\tau'$ be a trivalent graph on $\Sigma_g$ dual to a triangulation $\tau$ of $\Sigma_g$ with one vertex. Let $e_1',e_2',\cdots,e_{6g-3}'$ be a numbering of the edges of $\tau'$ such that the last three edges are associated to a  non-separating vertex and let $x_1, \cdots,x_{6g-3}$ be the cross ratios of the edges.  The projection of the cross ratio parameter space ${\mathcal C}_{\tau}$ to ${\mathbb R}^{6g-6}$  given by taking the first $6g-6$ variables,  
$$p: {\mathcal C}_{\tau} \longrightarrow {\mathbb R}^{6g-6}, $$
$$p(x_1,\cdots,x_{6g-3})=(x_1,\cdots,x_{6g-6})$$
is a homeomorphism onto its image. Furthermore, the  image of the projection is a strictly convex subset of ${\mathbb R}^{6g-6}$ with the property that if $(a_1,\cdots,a_{6g-6}) \in p({\mathcal C}_{\tau})$, then $(b_1, \cdots,b_{6g-6}) \in p({\mathcal C}_{\tau})$ whenever $b_i \geq a_i$ for $i=1,\cdots,6g-6$.
\end{lemma}

\medskip

\noindent {\it Proof:} We adopt the notation of proposition \ref{Prop:Derivative} and lemma \ref{Lem:Manifold}, and denote the last $3$ variables by $x$, $y$ and $z$ respectively. 

If a vector  ${\bf x} = (x_1,\cdots, x_{6g-6},x,y,z)$  represents a point in ${\mathcal C}_{\tau}$, the associated word is 
{\scriptsize $$W_{\bf x}={\left(\begin{array}{cc} 0 & 1 \\ -1 & x\\ \end{array}\right)} 
T
{\left(\begin{array}{cc} 0 & 1 \\ -1 & x\\ \end{array}\right)}
{\left(\begin{array}{cc} 0 & 1 \\ -1 & y\\ \end{array}\right)}
U
{\left(\begin{array}{cc} 0 & 1 \\ -1 & y\\ \end{array}\right)}
{\left(\begin{array}{cc} 0 & 1 \\ -1 & z\\ \end{array}\right)}
V
{\left(\begin{array}{cc} 0 & 1 \\ -1 & z\\ \end{array}\right)}$$}%
where $T$, $U$ and $V$ are subwords of $W$ of lengths $j-2$, $i-j-3$ and $12g-i-7$ respectively that depend only on the parameters $x_1,\cdots,x_{6g-6}$. By the main lemma (lemma \ref{Lem:MainLemma}), $T$, $U$ and $V$ are strictly admissible. \\
Now suppose that $(x_1,\cdots, x_{6g-6})$ are chosen so that $T$, $U$ and $V$ are strictly admissible. 

\smallskip

\noindent {\bf Claim:}
{\it  There exists a unique triple $(x,y,z) \in {\mathbb R}_+^3$ such that $W_{\bf x}$ satisfies the conditions of lemma \ref{Lem:MainLemma}, that is, $W_{\bf x}=-I$ and all proper subwords of  $W_{\bf x}$  are admissible. }

\medskip

\noindent {\it Proof of claim:} Let $T=\!{\left(\begin{array}{cc} t_1 & t_2 \\ t_3 & t_4\\ \end{array}\right)}$, $U=\!{\left(\begin{array}{cc} u_1 & u_2 \\ u_3 & u_4\\ \end{array}\right)}$, $V=\!{\left(\begin{array}{cc} v_1 & v_2 \\ v_3 & v_4\\ \end{array}\right)}$, $X={\left(\begin{array}{cc} 0 & 1 \\ -1 & x\\ \end{array}\right)}$, $Y={\left(\begin{array}{cc} 0 & 1 \\ -1 & y\\ \end{array}\right)}$, $Z={\left(\begin{array}{cc} 0 & 1 \\ -1 & z\\ \end{array}\right)}$. Note that the entries of each of $T$, $U$ and $V$ depend only on $x_1,\cdots,x_{6g-6}$ and that $t_1,u_1,v_1 \leq 0$, $t_2,u_2,v_2 >0$, $t_3,u_3,v_3 <0$ and $t_4,u_4,v_4>0$ since $T$, $U$ and $V$ are strictly admissible. 

By proposition \ref{Prop:ByX} (iii), $XTX$, $YUY$ and $ZVZ$ are strictly admissible if and only if $x>\alpha$, $y>\beta$, and $z>\gamma$ where 
$${\displaystyle \alpha= \frac{t_2-t_3+\sqrt{(t_2+t_3)^2+4}}{2t_4} , \qquad \beta= \frac{u_2-u_3+\sqrt{(u_2+u_3)^2+4}}{2u_4},}$$ $${\displaystyle \gamma= \frac{v_2-v_3+\sqrt{(v_2+v_3)^2+4}}{2v_4}}.$$
Since we require all subwords of $W_{\bf x}$  to be  admissible, we are only interested in solutions of $W_{\bf x}=-I$ where $x >\alpha$, $y >\beta$ and $z>\gamma$. We will show that there exists a unique triple $(x,y,z)$ satisfying $x >\alpha$, $y >\beta$ and $z>\gamma$ such that 
\begin{eqnarray*}
 \textrm{(2,2) entries of }
\left\{
 \begin{array}{l}
  XTXYUY=-(ZVZ)^{-1}
  \\
  YUYZVZ=-(XTX)^{-1}
  \\
  ZVZXTX=-(YUY)^{-1}
 \end{array}
\right.
.
\end{eqnarray*}

\noindent We start with the first equation. \!Equating the $(2,2)$ entry of $XTXYUY$ with that of $-(ZVZ)^{-1}$, we get
{\small $$ (t_2-t_4x)(u_3+u_4y)+(t_4x^2+(t_3-t_2)x-t_1)(u_4y^2+(u_3-u_2)y-u_1)  =v_4. \eqno(*)$$}
Let $S_1$ be the surface in ${\mathbb R}^3$
defined by $(*)$ and the inequalities $x>\alpha$, $y>\beta$, $z>\gamma$. Since the equation is independent of $z$, the surface is parallel to the $z$-axis and it suffices to analyze the intersection of the surface with any horizontal plane $z=z_0>\gamma$. Denote the left hand side of $(*)$ by $h(x,y)$ and note that it represents the $(2,2)$ entry of $XTXYUY$ and is quadratic in both $x$ and $y$. 

Fix $y=y_0>\beta$. Then $h(x,y_0)$ is a quadratic polynomial in $x$ with positive leading coefficient so that $\lim_{x \rightarrow \infty}h(x,y_0)=\infty$. At $x=\alpha$, 
\begin{eqnarray*}
XTXYUY&=&{\left(\begin{array}{cc} - & + \\ -& 0\\ \end{array}\right)}
{\left(\begin{array}{cc} -u_4 & u_3+u_4y_0 \\ u_2-u_4y_0 & -u_1+(u_3-u_2)y_0+u_4y_0^2\\ \end{array}\right)}\\
&=&{\left(\begin{array}{cc} * & * \\ *& -\\ \end{array}\right)} \qquad \hbox{since} \qquad u_3+u_4y_0>0
\end{eqnarray*}
so $h(\alpha,y_0)<0$.
It follows that for each $y_0>\beta$, there exists a unique $x>\alpha$ such that $h(x,y_0)=v_4$ since $v_4>0$ (the other solution must be $<\alpha$).

In fact, more is true. For this particular solution, $XTXYUY$ is strictly admissible. We first show that for $y_0>\beta$ and $x$ sufficiently large, $XTXYUY$ is strictly admissible. It is convenient to expand $XTXYUY$ as a word in the associated matrices  $A_1, A_2 \cdots$ where $A_1=X=A_j$, $A_{j+1}=Y=A_{i-1}$, $T=A_{n_2}\cdots A_{n_{j-1}}$, $U=A_{n_{j+2}}\cdots A_{n_{i-2}}$.  Then look at the circle configuration $\{ \overline C, C_0, \cdots , C_i\}$ corresponding to $XTXYUY$ (cf. definition \ref{Def:Admissibility}) with a suitable normalization. We choose the normalization so that the central circle $\overline{C}$  is the real line, $C_j$ is the line $\{Im(z)=1\}$ and $C_{j+1}$ is the circle $|z-\sqrt{-1}/2| = 1/2$ so that the interstice enclosed by $\overline{C}$, $C_j$ and $C_{j+1}$  is the standard interstice (see Figure \ref{Fig:Normalization}). Denote by $p_n$ the tangency between $\overline{C}$ and $C_n$, note that $p_n \in {\mathbb R}\cup \{ \infty \}$. We have:
\begin{itemize}
\item[(i)]
	$p_j=\infty$, $p_{j+1}=0$, by our normalization;
\item[(ii)] 
	$p_{j-1}=x$, since $A_j=X$; 
\item[(iii)]
	$0=p_{j+1}<p_{j+2}< \cdots <p_i<\infty $, since  $YUY=A_{j+1} \cdots A_{i-1}$ is strictly admissible;
\item[(iv)]
	$0<p_0<p_1<\cdots <p_{j-1}=x$, since $XTX=A_1 \cdots A_{j}$ is strictly admissible.
\end{itemize} 

\begin{figure}[ht]
  \begin{center}
    \begin{picture}(280,190)
      \put(0,120){\scalebox{0.7}{\includegraphics{Normalization-1.eps}}}
      \put(0,25){\scalebox{0.7}{\includegraphics{Normalization-2.eps}}}
\put(15,125){$p_{j+1}$}
\put(40,125){$p_{j+2}$}
\put(73,123){$\cdots$}
\put(105,125){$p_i$}
\put(150,125){$p_0$}
\put(177,123){$\cdots$}
\put(205,125){$p_{j-2}$}
\put(233,125){$p_{j-1}$}
\put(2,145){${\mathcal{I}_s}$}
\put(16,145){$C_{j+1}$}
\put(65,173){\vector(-1,-2){15}}
\put(63,176){$C_{j+2}$}
\put(122,169){\vector(-1,-2){15}}
\put(122,172){$C_{i}$}
\put(152,168){\vector(0,-1){30}}
\put(148,172){$C_{0}$}
\put(195,173){\vector(1,-2){15}}
\put(186,177){$C_{j-2}$}
\put(223,148){$C_{j-1}$}
\put(-11,133){$\overline{C}$}
\put(-15,166){$C_j$}
\put(50,106){(a) $XTXYUY$ is strictly admissible.}
\put(15,30){$p_{j+1}$}
\put(40,30){$p_{j+2}$}
\put(70,28){$\cdots$}
\put(109,39){\line(1,-2){6}}
\put(113,20){$p_i$}
\put(102.5,39){\line(-1,-2){6}}
\put(91,20){$p_0$}
\put(128,28){$\cdots$}
\put(152,30){$p_{j-2}$}
\put(175,30){$p_{j-1}$}
\put(2,50){${\mathcal{I}_s}$}
\put(16,50){$C_{j+1}$}
\put(65,78){\vector(-1,-2){15}}
\put(63,81){$C_{j+2}$}
\put(146,78){\vector(1,-2){15}}
\put(137,83){$C_{j-2}$}
\put(170,51){$C_{j-1}$}
\put(-11,38){$\overline{C}$}
\put(-15,71){$C_j$}
\put(50,-2){(b) $XTXYUY$ is not admissible.}
    \end{picture}
  \end{center}
  \caption[]{}
\label{Fig:Normalization}
\end{figure}

The condition that $XTXYUY$ is strictly admissible is the same as saying  the interstices do not intersect along $\overline{C}$.  This is easily seen to be equivalent to the condition that $p_i<p_0$. Now keeping all other cross ratios fixed, increasing $x$ by some positive number $c$ has the following effect on the points $\{ p_{k} \}$: 
\begin{itemize} 
\item[(a)]  
	$p_{j}, p_{j+1}, \cdots,p_{i}$ are left unchanged;
\item[(b)] 
	$p_1,p_2,\cdots,p_{j-1}$ are all shifted to the right by $c$;
\item[(c)]
	$p_0$ is shifted to the right by more than $c$ since the radius of the circle $C_0$ is decreased. 
\end{itemize}
It follows that for $x$ sufficiently large and $y>\beta$, $XTXYUY$ is strictly admissible. Furthermore, the above shows that there is a smallest value $x_1$ of $x$ when the word first becomes admissible but not strictly admissible, this occurs when $p_0=p_i$. In this case, $h(x_1,y_0)=0$. We have $x_1>\alpha$ since $XTX$ as a proper subword of an admissible word is strictly admissible. 

Since $h(x,y_0)$ is continuous and approaches infinity when $x \longrightarrow \infty$, it follows that there exists $x_2>x_1>\alpha$ for which $h(x_2,y_0)=v_4$, this must be the value we  have already found since it was the unique solution with $x>\alpha$. Since $x_2>x_1$, $XTXYUY$ is strictly admissible.

\smallskip

We now use the strict admissibility of $XTXYUY$ to show that the curve defined by $h(x,y)=v_4$ with $x>\alpha$ and $y>\beta$  is the graph of a decreasing function with vertical asymptote $x=\alpha$ and horizontal asymptote $y=\beta$. The gradient $\nabla h=\langle h_x,h_y \rangle$ at a point on the curve $h(x,y)=v_4$ can be calculated using the product rule as in the proof of proposition \ref{Prop:Derivative}. We have $h_x$ given by the $(2,2)$ entry of the matrix
\begin{eqnarray*}
&&{\left(\begin{array}{cc} 0 & 0 \\ 0& 1\\ \end{array}\right)}TXYUY+XT{\left(\begin{array}{cc} 0 & 0 \\ 0& 1\\ \end{array}\right)}YUY \\
&=&
{\left(\begin{array}{cc} 0 & 0 \\ 0& 1\\ \end{array}\right)}
{\left(\begin{array}{cc} - & + \\ -& +\\ \end{array}\right)}
+{\left(\begin{array}{cc} - &  +\\ -& +\\ \end{array}\right)}
{\left(\begin{array}{cc} 0 & 0 \\ 0& 1\\ \end{array}\right)}
{\left(\begin{array}{cc} - & + \\ -& +\\ \end{array}\right)}\\
&=&
{\left(\begin{array}{cc} 0 & 0 \\ -& +\\ \end{array}\right)}
+{\left(\begin{array}{cc} - &  +\\ -& +\\ \end{array}\right)}={\left(\begin{array}{cc} - &  +\\ -& +\\ \end{array}\right)}
\\
\end{eqnarray*}
  where the $(+)$ and $(-)$  entries follows from the fact that $XTXYUY$ is strictly admissible and hence so are all subwords. Similarly, $h_y$ is given by the $(2,2)$ entry of the matrix
$$XTX{\left(\begin{array}{cc} 0 & 0 \\ 0& 1\\ \end{array}\right)}UY+XTXYU{\left(\begin{array}{cc} 0 & 0 \\ 0& 1\\ \end{array}\right)} =
{\left(\begin{array}{cc} - &  +\\ -& +\\ \end{array}\right)}
$$
Hence, $\nabla h=\langle h_x,h_y \rangle=\langle +,+\rangle$ so  $h(x,y)=v_4$ defines $y$ as a decreasing function of $x$ for $x>\alpha$, $y>\beta$. The fact that $y=\beta$ is a horizontal asymptote and $x=\alpha$ is a vertical asymptote follows from an easy analysis of $h(x,y)$ as $x \longrightarrow \alpha^+$ and $y \longrightarrow \beta^+$.

\medskip
Clearly, the same analysis can be applied to the other two surfaces $S_2$ and $S_3$ given by $x>\alpha$, $y>\beta$, $z>\gamma$ and the equations
\begin{eqnarray*}
 \textrm{(2,2) entries of }
\left\{
 \begin{array}{l}
    YUYZVZ=-(XTX)^{-1}
  \\
  ZVZXTX=-(YUY)^{-1}
 \end{array}
\right.
\end{eqnarray*}
respectively. 
$S_2$ is given by an equation of the form $h_2(y,z)=t_4$ which defines $z$ as a decreasing function of $y$ and $S_3$ is given by an equation of the form $h_3(z,x)=u_4$ which defines $x$ as a decreasing function of $z$.
The intersection of $S_2$ and $S_3$ is a curve parameterized by $z \in (\gamma, \infty)$ whose projection onto the $xy$-plane is a curve with positive slope which approaches the point $(\alpha, \beta)$ as $z \longrightarrow \infty$ and whose $x$ and $y$ coordinates approach $\infty$ as $z \longrightarrow \gamma^+$. Since $h(x,y)=v_4$ has negative slope with asymptotes at $x=\alpha$ and $y=\beta$, there is a unique point of intersection of $S_1$ with the curve of intersection of $S_2$ and $S_3$, hence of the 3 surfaces. To complete the proof of the claim, we apply proposition \ref{Prop:MatrixIdentity} with $A=XTX$, $B=YUY$ and $C=ZVZ$.  
\qed 
\smallskip

 The uniqueness part of the claim now implies that the projection map $p$ is one-to-one. The existence part implies that in fact, ${\mathcal C}_{\tau}$ is homeomorphic to
$$\{(x_1,\cdots,x_{6g-6}) \in {\mathbb R}_+^{6g-6} ~|~ T,U, V \; \hbox{are strictly admissible} \}.$$ Now the last part of the lemma follows from the convexity properties proven in \ref{Prop:Convex}.
\qed

%
%

\section{Rigidity for one circle packings}


\subsection{Rigidity}

In the previous section, we gave a complete description of the cross ratio parameter space ${\mathcal C}_{\tau}$ for a triangulation $\tau$ of $\Sigma_g$ with one vertex. By definition, each element of ${\mathcal C}_{\tau}$ determines a projective structure on $\Sigma_g$ together with a circle packing with nerve $\tau$. The forgetting map $f : {\mathcal C}_{\tau} \longrightarrow {\mathcal P}_g$  is an embedding if we can show that the circle packings are rigid, that is, given any two circle packings $P_1$ and $P_2$ on a surface with projective structure such that their nerves are isotopic, then there exists a projective automorphism of the surface isotopic to the identity taking $P_1$ to $P_2$. 

There have been many recent papers concerning rigidity results of circle packings, the most general perhaps being \cite{He} which also contains a comprehensive list of references.
 However, the results do not seem to apply to our problem, indeed, many of the methods do not generalize and many basic results like the ring lemma do not apply in our case. We prove the following result which implies rigidity of circle packings in the one-circle case. The main theorem then follows as a corollary of lemma \ref{Lem:Cell} and lemma \ref{Lem:Rigid}.

\medskip

\begin{lemma}\label{Lem:Rigid}
Suppose that $\tau$ is a triangulation of $\Sigma_g$ with one vertex and ${\bf c}, {\bf c'} \in {\mathcal C}_{\tau}$. Then if $f({\bf c}) = f({\bf c'})$, in other words,  ${\bf c}$  and ${\bf c'}$ define the same projective structure on $\Sigma_g$, then ${\bf c}={\bf c'}$. Hence, the circle packings of $\Sigma_g$ by one circle are rigid.
\end{lemma}

\medskip
\noindent {\it Proof:} We adopt the notation of proposition \ref{Prop:Derivative}, lemma \ref{Lem:Manifold} and \ref{Lem:Cell} so that the  $3$  values of  ${\bf c}$ (respectively ${\bf c'}$) are denoted by $x$, $y$, $z$  (respectively $x'$, $y'$, $z'$)  and are the cross ratios of three edges of the dual graph $\tau'$ meeting at a non-separating vertex $v'$. 
The corresponding words of associated matrices for ${\bf c}$ and ${\bf c'}$ are given by
$$
W= XTXYUYZVZ, \qquad W'=X'T'X'Y'U'Y'Z'V'Z'
$$
where $X={\left(\begin{array}{cc} 0 & 1 \\ -1& x\\ \end{array}\right)}$,
$Y={\left(\begin{array}{cc} 0 & 1 \\ -1& y\\ \end{array}\right)}$,
$Z={\left(\begin{array}{cc} 0 & 1 \\ -1& z\\ \end{array}\right)}$, and 
$T={\left(\begin{array}{cc} t_1 & t_2 \\ t_3& t_4\\ \end{array}\right)}$, $U={\left(\begin{array}{cc} u_1 & u_2 \\ u_3& u_4\\ \end{array}\right)}$ and $V={\left(\begin{array}{cc} v_1 & v_2 \\ v_3& v_4\\ \end{array}\right)}$ are subwords of $W$ of lengths $j-2$, $i-j-3$ and $12g-i-7$ respectively and depend only on $6g-6$ variables, and similarly for $X'$, $Y'$, $Z'$, $T'$, $U'$ and $V'$. We note that these matrices are in $SL_2({\mathbb C})$ but when we consider the holonomy representation, we consider them also as elements in $PSL_2({\mathbb C})$, there should be no confusion as the choice should be clear from the  context. 

We first fix the developing map and hence the holonomy representation for both ${\bf c}$ and ${\bf c'}$ by developing the circle to the real line and the interstice containing $v'$ to the standard interstice. We get a configuration of circles $\{\overline{C}, C_1, \cdots,C_{12g-6}\}$ where $\overline{C}$ is the real line, $C_1$ is the circle $|z-\sqrt{-1}/2|=1/2$ and $C_{12g-6}$ is the line $Im(z)=1$ for ${\bf c}$ and a similar configuration for ${\bf c'}$. Let $\gamma_1, \gamma_2, \gamma_3 \in \pi_1(\Sigma_g)$ where  $\gamma_1:s_1 \rightarrow s_j$, $\gamma_2:s_{j+1}\rightarrow s_{i-1}$, $\gamma_3:s_{i}\rightarrow s_{12g-6}$.
Then $\gamma_3 \gamma_2 \gamma_1=id$ and if $\rho$ and $\rho'$ are the holonomy representations corresponding to ${\bf c}$  and ${\bf c'}$ respectively, then by proposition \ref{Prop:Holonomy},
\begin{eqnarray*}
\rho(\gamma_1)=XTXR^{-1}, & \quad \rho(\gamma_2)=XTXYUYR^{-1}(XTX)^{-1}, \\
 \rho'(\gamma_1)=X'T'X'R^{-1}, & \quad \rho'(\gamma_2)=X'T'X'Y'U'Y'R^{-1}(X'T'X')^{-1}
\end{eqnarray*}
where $R^{-1}={\left(\begin{array}{cc} 1 & -\sqrt{-1} \\ -\sqrt{-1} & 0\\ \end{array}\right)}$.
Now  since ${\bf c}$ and ${\bf c'}$ give the same projective structures, the holonomy representations $\rho$ and $\rho'$ must be conjugate, that is, there exists some $H \in PSL_2({\mathbb C})$ such that for any $\gamma \in \pi_1(\Sigma_\gamma)$, $\rho'(\gamma)=H^{-1}\rho(\gamma)H$. Hence
\begin{eqnarray*}
tr(\rho(\gamma_1))&=&\pm tr(\rho'(\gamma_1))\\
\Rightarrow tr(XTXR^{-1})&=&\pm tr(X'T'X'R^{-1})\\
\Rightarrow    -t_4 -(t_2+t_3)\sqrt{-1}         &=&\pm (-t_4' -(t_2'+t_3')\sqrt{-1}).
\end{eqnarray*}
Equating real and imaginary parts, and using the fact that $t_4, t_4'>0$, we get
$t_4=t_4'$ and $t_2+t_3=t_2'+t_3'$. Note that since $t_4>0$, $tr(\rho(\gamma_1)) \neq 0$ and so $\rho(\gamma_1)$ is not an elliptic element of order 2. Let $t_2'=t_2-k$, so $t_3'=t_3+k$.
Then
\begin{align}
X'T'X'&= {\left(\begin{array}{cc} 0 & 1 \\ -1& x'\\ \end{array}\right)}
{\left(\begin{array}{cc} t_1' & t_2' \\ t_3'& t_4'\\ \end{array}\right)}
{\left(\begin{array}{cc} 0 & 1 \\ -1& x'\\ \end{array}\right)} \displaybreak[0] \notag \\
&={\left(\begin{array}{cc} -t_4' & t_3'+t_4'x' \\ t_2'-t_4'x' & t_4'{x'}^2+(t_3'-t_2')x'-t_1'\\ \end{array}\right)} \displaybreak[0] \notag \\
&={\left(\begin{array}{cc} -t_4 & t_3+(t_4x'+k) \\ t_2-(t_4x'+k) & t_4{x'}^2+(t_3-t_2+2k)x'-t_1'\\ \end{array}\right)} \displaybreak[0] \notag \\
&= {\left(\begin{array}{cc} 0 & 1 \\ -1& x''\\ \end{array}\right)}
{\left(\begin{array}{cc} t_1 & t_2 \\ t_3& t_4\\ \end{array}\right)}
{\left(\begin{array}{cc} 0 & 1 \\ -1& x''\\ \end{array}\right)} \displaybreak[0] \notag \\
&=X''TX''\qquad \hbox{where} \quad x''=x'+\frac{k}{t_4}, \quad X''={\left(\begin{array}{cc} 0 & 1 \\ -1& x''\\ \end{array}\right).} \displaybreak[0] \notag
\end{align}
$X''TX''$ is strictly admissible since $X'T'X'$ is strictly admissible (admissibility of $X''T$ and $TX''$ follows from proposition \ref{Prop:ByX} by looking at the $(1,2)$ entry and the $(2,1)$ entries and admissibility of $X''TX''$ then follows from the $(2,2)$ entry). Similarly, we have 
\begin{eqnarray*}
Y'U'Y'&=& {\left(\begin{array}{cc} 0 & 1 \\ -1& y''\\ \end{array}\right)}
{\left(\begin{array}{cc} u_1 & u_2 \\ u_3& u_4\\ \end{array}\right)}
{\left(\begin{array}{cc} 0 & 1 \\ -1& y''\\ \end{array}\right)} =Y''UY''\\
\hbox{and}\qquad Z'V'Z'&=& {\left(\begin{array}{cc} 0 & 1 \\ -1& z''\\ \end{array}\right)}
{\left(\begin{array}{cc} v_1 & v_2 \\ v_3& v_4\\ \end{array}\right)}
{\left(\begin{array}{cc} 0 & 1 \\ -1& z''\\ \end{array}\right)} =Z''VZ''.\\
\end{eqnarray*}

By the uniqueness part of the claim in the proof of lemma \ref{Lem:Cell}, we get that $x''=x$, $y''=y$ and $z''=z$ so that $XTX=X'T'X'$, $YUY=Y'U'Y'$ and $ZVZ=Z'V'Z'$. Hence $\rho(\gamma_1)=\rho'(\gamma_1)$, $\rho(\gamma_2)=\rho'(\gamma_2)$ and $\rho(\gamma_3)=\rho'(\gamma_3)$, that is, the conjugating element $H$ commutes with $\rho(\gamma_1), \rho(\gamma_2)$ and $\rho(\gamma_3)$.

We now claim that $\rho(\gamma_1)$ and $\rho(\gamma_2)$ do not commute. If true, then $H=I$ since it commutes with both $\rho(\gamma_1)$ and $\rho(\gamma_2)$, and  they are  not elliptics of order 2.  Hence $\rho=\rho'$ and the two configurations for ${\bf c}$ and ${\bf c'}$ are identical and so ${\bf c}={\bf c'}$.

It remains to prove that $\rho(\gamma_1)$ and $\rho(\gamma_2)$ do not commute (equivalently,  $\rho(\gamma_1)$, $\rho(\gamma_2)$ and $\rho(\gamma_3)$ do not pairwise commute since $\gamma_3 \gamma_2 \gamma_1=id$). \\
Suppose not. Consider the configuration of circles $\{\overline{C}, C_1, \cdots, C_{12g-6}\}$. Recall that the tangency points $\{p_m\}$ of $\overline{C}$ with $C_m$  satisfy 
$$0=p_1<p_j<p_{j+1}<p_{i-1}<p_i<p_{12g-6}=\infty$$
 with $\rho(\gamma_1)(p_1)=p_j$, $\rho(\gamma_2)(p_{j+1})=p_{i-1}$, $\rho(\gamma_3)(p_{i})=p_{12g-6}$. \\
 We have
$\rho(\gamma_2)(\overline{C})=C_{i-1}$ is tangent to $\overline{C}$, hence 
$\rho(\gamma_1)(\rho(\gamma_2)(\overline{C}))$ is tangent to $\rho(\gamma_1)(\overline{C})=C_j$.
By assumption $\rho(\gamma_1)\rho(\gamma_2)=\rho(\gamma_2)\rho(\gamma_1)$ so $$\rho(\gamma_1)(\rho(\gamma_2)(\overline{C}))=\rho(\gamma_2)\rho(\gamma_1)(\overline{C})=\rho(\gamma_3^{-1})(\overline{C})=C_i$$ is tangent to $C_j$. Similarly, permuting the roles of $\gamma_1$, $\gamma_2$ and $\gamma_3$, the commuting condition implies that $C_{j+1}$ is tangent to $C_{12g-6}$ and $C_{i-1}$ is tangent to $C_1$. Now $C_{j+1}$ is tangent to $C_{12g-6}$ and $\overline{C}$ implies that $C_{j+1}=\{z \in {\mathbb C}: |z-(a+\sqrt{-1}/2)|=1/2\}$, the circle of radius $1/2$ with center at $a+\sqrt{-1}/2$ where $p_{j+1}=a>0$ (see Figure \ref{Fig:Tangent}). Let $r$ be the radius of the circle in the upper half plane tangent to the real line at the origin and also tangent to $C_{j+1}$. Since $C_j$ is tangent to $C_{j+1}$ and its point of tangency with the real line lies between $0$ and $p_{j+1}$, the radius of $C_j$ is less than $r$. Now since $C_{i-1}$ is tangent to $C_1$ and its point of tangency $p_{i-1}$ with the real line lies between $p_{j+1}$ and $\infty$, the radius of  $C_{i-1}$ is greater than $r$. Hence the radius of  $C_j$  is strictly less than the radius of  $C_{i-1}$  with the points of tangencies satisfying $p_j<p_{i-1}$. Now $C_i$ is tangent to $C_{i-1}$ and its tangency point with the real line lies between $p_{i-1}$ and $\infty$ which implies that it does not intersect $p_{j}$, giving a contradiction. Hence $\rho(\gamma_1)$ and $\rho(\gamma_2)$ do not commute. \qed

\begin{figure}[ht]
  \begin{center}
    \begin{picture}(280,70)
      \put(0,0){\scalebox{0.7}{\includegraphics{Tangent.eps}}}
\put(15,7){$p_{1}$}
\put(40,7){$p_{2}$}
\put(58,5){$\cdots$}
\put(79,7){$p_{j}$}
\put(96,7){$p_{j+1}$}
\put(120,7){$p_{j+2}$}
\put(144,5){$\cdots$}
\put(160,7){$p_{i-1}$}
\put(195,5){$\cdot\ \ \cdot\ \ \cdot\ \ $}
\put(250,7){$p_{i}$}
\put(22,25){$C_{1}$}
\put(65,54){\vector(-1,-2){15}}
\put(63,56){$C_{2}$}
\put(87,26){$C_{j+1}$}
\put(135,53){\vector(-1,-2){15}}
\put(133,55){$C_{j+2}$}
\put(211,25){$C_{i-1}$}
\put(267,50){\vector(-1,-2){15}}
\put(265,52){$C_{i}$}
\put(-11,13){$\overline{C}$}
\put(-35,43){$C_{12g-6}$}
    \end{picture}
  \end{center}
  \caption[]{}
\label{Fig:Tangent}
\end{figure}

\medskip


\subsection{Concluding remark}

\smallskip

Much of our arguments above rely on the combinatorics of the one circle packing and it remains unsolved whether the main theorem holds for a general triangulation $\tau$ of $\Sigma_g$. There are also other aspects to be explored, for example one can look at circle packings on surfaces with cone projective structures, where the cone points would be the centers of some of the circles of the packing. Some results have been obtained and details will be given in a future paper. Of particular interest is the case of branched circle packings, this correspond to cone points with cone angles which are multiples of $2\pi$. For example one can look at the deformation space of cone projective structures on the surface admitting a branched circle packing by one circle with branch point of multiplicity one, this would correspond to the condition that $W_v=I$ and that in the developing map, each circle is surrounded by $12g-6$ circles which wrap around the  circle exactly twice. It seems that in this case, the deformation space ${\mathcal C}_{\tau}$ may admit some interesting topology, analogous to the case of the components of the space of representations   $Hom( \pi_1(\Sigma_g), PSL_2({\mathbb R}))/PSL_2({\mathbb R})$   where the Euler class is not $\pm (2-2g)$. 

\bigskip


\section{Appendix}

We consider the case of the torus $\Sigma_1$ with a circle packing by one circle.  Here the nerve $\tau$ consists of one vertex $v$ and 3 edges $e_1$, $e_2$ and $e_3$ with cross ratios $x>0$, $y>0$ and $z>0$ respectively and associated matrices $X={\left(\begin{array}{cc} 0 & 1 \\ -1& x\\ \end{array}\right)}$, $Y={\left(\begin{array}{cc} 0 & 1 \\ -1& y\\ \end{array}\right)}$ and $Z={\left(\begin{array}{cc} 0 & 1 \\ -1& z\\ \end{array}\right)}$. The word associated to the vertex is given by $W=XYZXYZ$ and by condition $(1)$ of lemma \ref{Lem:MainLemma},  
\begin{eqnarray*} 
XYZXYZ&=&-I\\ 
\Leftrightarrow XYZ&=&-(XYZ)^{-1}\\
\Leftrightarrow {\left(\begin{array}{cc} -y & yz-1 \\ 1-xy& xyz-x-z\\ \end{array}\right)}&=&{\left(\begin{array}{cc} x+z-xyz & yz-1 \\ 1-xy & y\\ \end{array}\right)}\\
\Leftrightarrow \qquad xyz &=&x+y+z. \qquad \qquad \qquad
\qquad (**)
\end{eqnarray*}
Note that three equations derived from the matrix identity reduces to just one equation in this case.
 Condition $(2)$ of lemma \ref{Lem:MainLemma} is automatically satisfied if $x$, $y$ and $z$ are all positive and satisfy $(**)$ since $xy-1=(x+y)/z>0$, $yz-1=(y+z)/x>0$,
$zx-1=(z+x)/y>0$,
and $xyz-x-z=y>0$ hence $XYZ$ is strictly admissible and so $XYZXYZ$ corresponds to a configuration where the surrounding circles wrap around the central circle exactly once.
The cross ratio parameter space ${\mathcal C}_{\tau}$ is hence given by 
$${\mathcal C}_{\tau}=\{(x,y,z) \in {\mathbb R}^3 |~xyz=x+y+z, \, x,y,z > 0 \}$$
and the image of the projection to the $xy$-plane is the set 
$$p({\mathcal C}_{\tau})=\{(x,y) \in {\mathbb R}^2 |~xy-1>0, \, x,y > 0 \}$$
which is convex.

To show that the circle packings are rigid, we directly compute the trace of the holonomy representation of a  pair of generating elements of $\pi_1(\Sigma_1)$ (note that $\pi_1(\Sigma_1)$ is abelian and so the holonomy image of $\pi_1(\Sigma_1)$ in $PSL_2({\mathbb C})$ is necessarily an elementary group unlike the higher genus case).
From example \ref{Ex:Holonomy}, $\pi_1(\Sigma_1)$ is generated by $\gamma_1$ and $\gamma_2$ with
\begin{eqnarray*}
\rho(\gamma_1)&=& ZXR ={\left(\begin{array}{cc} x\sqrt{-1} & x-\sqrt{-1} \\ (xz-1)\sqrt{-1}& (xz-1)-z\sqrt{-1}\\ \end{array}\right)}\\
\rho(\gamma_2)&=& ZXYRZ^{-1} =Z{\left(\begin{array}{cc} y\sqrt{-1} & y-\sqrt{-1} \\ (yx-1)\sqrt{-1}& (yx-1)-x\sqrt{-1}\\ \end{array}\right)}Z^{-1} \\
\end{eqnarray*}
so $\hbox{tr}(\rho(\gamma_1))=\pm (xz-1+(x-z)\sqrt{-1})$, \quad $\hbox{tr}(\rho(\gamma_2))=\pm(yx-1+(y-x)\sqrt{-1})$.
If $(x',y',z')$ is another point in the cross ratio parameter space which gives the same projective structure, then equating the traces, we get
\begin{eqnarray*}
x'z'-1+(x'-z')\sqrt{-1}&=&\pm (xz-1+(x-z)\sqrt{-1}),\\
 \quad y'x'-1+(y'-x')\sqrt{-1}&=&\pm(yx-1+(y-x)\sqrt{-1}). 
\end{eqnarray*}
Equating real and imaginary parts and taking into account the inequalities satisfied by the variables, we get $x'=x$, $y'=y$ and $z'=z$, hence the circle packing is rigid.   
It follows that the space of projective structures on the torus admitting a circle packing by one circle is homeomorphic to ${\mathbb R}^2$. Finally, we note that a projective structure on the torus is necessarily an euclidean or similarity structure, in the case of the euclidean structure, it is the one with fundamental domain a regular hexagon. In the other cases, it can be easily shown that the developing maps of the similarity structures obtained in this way are precisely the well-known Doyle spirals.

%
%

\end{document}